\documentclass[12pt]{article}

\usepackage[cp1251]{inputenc}
\usepackage[T2A]{fontenc}
\usepackage[english]{babel}

\usepackage{amsfonts,amsmath,amssymb, amsthm}
\usepackage{graphicx}

\usepackage{amsmath,amsthm,amssymb,graphicx}
\usepackage{rotating}
\usepackage{epsfig}
\usepackage{lscape}
\usepackage{subfigure}

\newtheorem{theorem}{Theorem}
\newtheorem{cor}{Corollary}
\newtheorem{example}{Example}
\newtheorem{remark}{Remark}
\newcommand{\eps}{\varepsilon}
\newcommand{\R}{\mathbb R}
\newcommand{\T}{\mathbb T}

\newcommand{\dst}{\displaystyle}
\title{Double Hamiltonian Hopf Bifurcation: \\normalization
and normal form non-integrability}
\author{L.M. Lerman$^1$, R. Mazrooei-Sebdani$^{2, 3}$\footnote{This research was
in part supported by a grant from IPM (No.1403700317).}, N.E. Kulagin$^4$\\
\normalsize
$^1$ HSE University, Nizhny Novgorod Branch, Russia\\
\normalsize
$^2$ Department of Mathematical Sciences, Isfahan University of Technology, Iran,\\
\normalsize
$^3$ School of Mathematics, Institute for Research in Fundamental Sciences, Tehran, Iran\\
\normalsize
$^4$ A.N. Frumkin Institute of Physical Chemistry and Electrochemistry RAS, Moscow, Russia}
\date{}
 
\begin{document}
\maketitle

\begin{abstract}
The double Hamiltonian Hopf bifurcation is studied, i.e. a generic two-parametric unfolding
of a smooth Hamiltonian system with four degrees of freedom which has at the critical value
of parameters the equilibrium with two pairs of double non semi-simple pure imaginary
eigenvalues $\pm i\omega_1,$ $\pm i\omega_2,$ $\omega_1\ne \omega_2$ under an assumption
of absence of strong resonances between $\omega_1,\omega_2$. We derive the normal form of
the unfolding, when the ratio
$\omega_1/\omega_2$ is irrational and study the truncated normal form of the fourth
order. This truncated normal form is the same under the absence of strong resonances.
The normal form has two quadratic integrals generating a symplectic periodic
action of the abelian group $\T^2.$ After reduction by means of these integrals we come
to the reduced system with two degrees of freedom that is proven to be non-integrable
for almost all values of its coefficients. Integrable such systems are also possible at some
special values of coefficients, related examples are presented. Some investigations of this
truncated system are presented along with its bifurcations when varying small detuning
parameters. As an example of a system where this bifurcation is met, the system derived in
\cite{KuLe} is investigated. Its homoclinic solutions are examined numerically
when the system parameters correspond to a main equilibrium of the twofold saddle-focus
type.
\end{abstract}

\section{Introduction and set up}

The bifurcation, which name is in the title, has been met for the first
time (as we aware of) when studying patterns of the stationary Swift-Hohenberg equation
on the plane \cite{KuLe}. Using Bubnov-Galerkin approximation to search for periodically
modulated wall-type solutions, the system of 8 or 12 differential equations for modes was
derived. These systems are reversible Hamiltonian and their homoclinic solutions
of the equilibrium at the origin correspond to solutions were sought for. Such solutions
were found in \cite{KuLe} by means of some numerical simulations. Here we want to understand
some features of the bifurcation with some details. We are confident that this codimension
two bifurcation will be met in many applied systems, so it deserves a thorough
study. It would be very interesting to find this bifurcation in systems
from applications.

We intend to investigate the local dynamics and bifurcations near the equilibrium,
so the problem can be studied in the local symplectic coordinates. In order
to formulate the problem precisely, we consider the linear symplectic space
$(\R^8, \Omega)$ with its standard symplectic form $\Omega$ and suppose a smooth ($C^\infty$)
Hamilton function $H$ be defined in a neighborhood $U$ of the origin $O$ where $dH=0$.
This means the point $O$ be an equilibrium of the Hamiltonian vector field $X_H$. Let
$(x,y)$ be symplectic coordinates in the neighborhood $U$ of $O$,
$x = (x_1,\ldots,x_4),$ $y = (y_1,\ldots, y_4),$ $\Omega = \sum_k dy_k\wedge dx_k$, $H=H(x,y)$.
Then the vector field $X_H$ is written down in the standard form
$$
\dot x = H_y,\;\dot y = - H_x,
$$
and the linearized system at the equilibrium $O$ is as follows
$$
\dot \xi = H^0_{xy}\xi + H^0_{yy}\eta,\;\dot\eta = - H^0_{xx}\xi - H^0_{xy}\eta,
$$
where zeroth upper index means that the related derivatives are calculated at
the point $O$. Denote $J$ the standard $8\times 8$ matrix
$$
\begin{pmatrix}0&I_4\\-I_4&0\end{pmatrix}.
$$
Also denote $4\times 4$-matrices $A = H^0_{xx}$, $C = H^0_{yy},$ being
symmetric and $B = H^0_{xy},$ then the linearized system casts in the form
$$
\begin{pmatrix}\dot\xi\\\dot\eta\end{pmatrix} = JR\begin{pmatrix}\xi\\\eta\end{pmatrix},\;
R = \begin{pmatrix}A&B\\B^\top&C\end{pmatrix}.
$$

We suppose the Hamiltonian matrix $JR$ of the linear system have two pairs of double imaginary
eigenvalues $\pm i\omega_1, \pm i\omega_2$, $\omega_1\ne\omega_2,$ $\omega_i \ne 0,$ being
each non semi-simple. The operator $JR$ is decomposable:
the linear symplectic space $\R^8$ decomposes into direct sum of two invariant symplectic
four-dimensional subspaces in which the restriction of the operator $JR$ generates
indecomposable non semi-simple Hamiltonian operator that has two double non semi-simple
pure imaginary eigenvalues $\pm i\omega_1$ and $\pm i\omega_2$, respectively,
$\omega_1\ne \omega_2$.

Such equilibrium has a degeneracy of co-dimension two in the class of smooth Hamiltonians,
so it has to be examined in a generic two-parameter unfolding. We want to understand main
features of its local dynamics and study its bifurcations. The problem consists in the
interrelation of two Hamiltonian Hopf bifurcations \cite{Meer}, so we expect all phenomena
met there will come into play here. As a particular task, we want to investigate under which
conditions small homoclinic orbits of the equilibrium $O$ can arise in the unfolding. As is
known, the creation of small homoclinic orbits at the Hamiltonian Hopf bifurcation \cite{Meer}
occurs if some coefficient in the normal form of the fourth order be positive, whence it cannot
be decided via the linearization. The conditions of this phenomenon here are unknown, this
question was of interest for the system studied in \cite{KL1}. As was said, there such orbits
were found numerically.

It is worth noting from the very beginning that in the critical system the equilibrium $O$
has not zero eigenvalues, therefore this equilibrium persists and smoothly depends on
parameters, if the system depends smoothly on them. So, without loss of generality, one may
assume that for all systems of the unfolding the equilibrium $O$ does not depend on parameters
and is located at the origin. This will be assumed henceforth.

From our assumptions on the type of the equilibrium for the critical system, the usage of
results of \cite{Galin} (see also, \cite{Kocak,LM,Hovejin,Chow}) allows one to
transform the quadratic parts of the Hamiltonians of the unfolding at the equilibrium by
a linear symplectic change of variable smoothly depending on parameters to
the normal form with small parameters $\lambda_1, \lambda_2$
$$
H_2 = (p_2q_1 - \omega^2_1p_1q_2)\pm\frac12(\frac{q_1^2}{\omega_1^2} + q_2^2) +
(p_4q_3 - \omega^2_2p_3q_4)\pm \frac12(\frac{q_3^2}{\omega^2_2} + q_4^2) +
\frac{\lambda_1}{2}p_1^2 + \frac{\lambda_2}{2}p_3^2.
$$
We correct this form a bit by means of the additional linear symplectic change of variables:
$Q_1=q_1/\omega_1, P_1 = \omega_1p_1,$ $Q_3=q_3/\omega_2, P_3 = \omega_2p_3,$ after that
the Hamiltonians take  the form (we preserve the same notations for new variables)
$$
H_2 = \omega_1(p_2q_1 - p_1q_2)\pm\frac12(q_1^2 + q_2^2) +
\omega_2(p_4q_3 - p_3q_4)\pm \frac12(q_3^2 + q_4^2) +
\frac{\eps_1}{2}p_1^2 + \frac{\eps_2}{2}p_3^2,
$$
with new small parameters $\eps_1 = \lambda_1/\omega_1^2,$ $\eps_2 = \lambda_2/\omega_2^2.$

In fact, we alter the perturbation terms and actually consider $H_2$ in
the form
\begin{equation}\label{crit}
H_2 = \omega_1(p_2q_1 - p_1q_2)\pm\frac12(q_1^2 + q_2^2) +
\omega_2(p_4q_3 - p_3q_4)\pm \frac12(q_3^2 + q_4^2) +
\frac{\eps_1}{2}(p_1^2 + p_2^2) + \frac{\eps_2}{2}(p_3^2 + p_4^2),
\end{equation}
where coefficients $\omega_k, k=1,2,$ are smooth functions of parameters
$\eps_1,\eps_2.$ This and previous quadratic Hamiltonians can be
transformed to each other by a close to identity linear symplectic change
of variables and introducing new small parameters instead of $\eps_1,\eps_2.$

The critical (as $\eps_1=\eps_2=0$) quadratic Hamiltonian $H_2^0$ (\ref{crit}) can
be written as
\begin{equation}
\label{H2-2}
H_2^0 = H_2^s + H_2^{n},
\end{equation}
here
$$
H_2^s = \omega_1(p_2q_1 - p_1q_2)+\omega_2(p_4q_3 - p_3q_4),\;H_2^{n} =
\pm\frac12(q_1^2 + q_2^2) \pm \frac12(q_3^2 + q_4^2)
$$
represent the semi-simple part $H_2^s$ and its nilpotent part $H_2^n$.

The form (\ref{crit}) of the quadratic Hamiltonian shows that we deal with the
direct product of two double $1:-1$ resonance, i.e. $\omega_1:-\omega_1$ acting in
the invariant subspace $(q_1,q_2,p_1,p_2)$ and $\omega_2:-\omega_2$ acting in the invariant
subspace $(q_3,q_4,p_3,p_4)$. We also see that re-enumeration of variables allows one
to consider both numbers $\omega_1, \omega_2$ be positive and satisfying the inequality
$\omega_1 < \omega_2$, this will be assumed to hold in what follows. The passage of
the unfolding through a system that has an equilibrium with this resonance will be
called to be {\em double Hamiltonian Hopf bifurcation} or, briefly, DHHB.

Our study can be considered as a continuation of the rather long series of works performed
for the investigation of a more simple bifurcation, the so-called Hamiltonian Hopf bifurcation
named so after the title of the book \cite{Meer}. The first results there were obtained in
\cite{DH} and \cite{MSch} where the local families of periodic orbits were investigated.
The general problem of studying families of local periodic solutions near a resonant equilibrium
was put and partially solved in \cite{Duis}, though different partial cases were investigated
in many Russian papers, see, for instance, \cite{sok} and the book \cite{KhaSchnol}. These
approaches were extended onto systems with more degrees of freedom in \cite{HM}, see the
book \cite{Hanssmann}. The situation in non Hamiltonian systems is discussed in
\cite{BHW,GTtor}.

\section{Equilibrium types at the origin}\label{origin}

Before studying local bifurcations in the nonlinear system let us determine the types of
the equilibrium $O$ in the linear approximation for the
two-parametric unfolding. In the quadratic approximation, the Hamiltonian splits into two
independent partial quadratic Hamiltonians where the first acts in the invariant symplectic
subspace $q_3 = q_4 = p_3 = p_4 = 0$ with its coordinates $(q_1,q_2,p_1,p_2)$ and the second
one does in the invariant symplectic subspace $q_1 = q_2 = p_1 = p_2 = 0$ with its coordinates
$(q_3,q_4,p_3,p_4)$. Thus, the parameter $\eps_1$ governs by the Hamiltonian Hopf
Bifurcation for the first partial Hamiltonian and the parameter $\eps_2$ does the same for
the second partial Hamiltonian.
This shows that in the parameter plane $(\eps_1,\eps_2)$ within some sufficiently small disk
$P$ centered at the origin, we get four sectors being the intersection of $P$ with four
quadrants with all collections of signs for $(\eps_1,\eps_2)$. In either sector, the values
of signs are definite and we can calculate characteristic polynomials for both partial
Hamiltonians. If we omit the indices in $\eps_i$, then the characteristic equation
is as follows
$$
\lambda^4 + 2\lambda^2(\eps+\omega^2) + (\omega^2 - \eps)^2 = 0,
$$
where $\eps$, $\omega$ corresponds to $\eps_1$, $\omega_1$, when we consider the first partial
Hamiltonian, and to $\eps_2$, $\omega_2$, when we consider the second partial
Hamiltonian.

The sectors correspond to the following types of the equilibrium $O$:
\begin{itemize}
\item an elliptic point with eigenvalues $\pm [i(\omega_1\pm \sqrt{\eps_1}),
i(\omega_2\pm \sqrt{\eps_2})]$, when $\eps_1>0, \eps_2>0$;
\item an elliptic-saddle-focus type with eigenvalues $\pm [i(\omega_1\pm \sqrt{\eps_1}),
\pm\sqrt{-\eps_2}+ i\omega_2$)], when $\eps_1 > 0,$ $\eps_2 < 0$, and the same type
$\pm[\pm\sqrt{-\eps_1}+ i\omega_1, i(\omega_2\pm \sqrt{\eps_2})]$ for $\eps_1 < 0,
\eps_2 > 0$;
\item a twofold saddle-focus with eigenvalues $\pm[\pm\sqrt{-\eps_1}+ i\omega_1,
\pm\sqrt{-\eps_2}+ i\omega_2]$ as both $\eps_1, \eps_2$ negative.
\end{itemize}
The lines $\eps_1 = 0$ and $\eps_2 = 0$ are bifurcational, their crossing changes the type
of $O$.

Recall some elements of the local orbit structure for these structurally stable
equilibria within the class of smooth Hamiltonians. For the case of a twofold
saddle-focus, the equilibrium is of the hyperbolic type, it has locally two smooth invariant
four-dimensional manifolds, stable $W^s$ and unstable $W^u$ ones, containing all orbits
which tend to $O$ as $t\to \infty$ (for $W^s$) and $t\to -\infty$ (for $W^u$).
All other orbits in some sufficiently small (of the order $\sqrt{\eps_1^2+\eps_2^2}$)
neighborhood $U$ of $O$ with initial points outside of the union $W^s\cup W^u$ leave $U$
in both directions of time.
Here one of natural questions is as follows: both these local $W^s, W^u$ belong to the level
$H=H(O)$ being a smooth 7-dimensional submanifold everywhere except for the point $O$ where it
has a conical singularity, can these submanifolds intersect each other (transversely generically)
under continuation by the flow staying nevertheless in a small neighborhood of $O$ (but larger
than $U$, say of the order $\sqrt{|\eps_1|+|\eps_2|}$), as this occurs for the Hamiltonian
Hopf bifurcation \cite{Meer}, when the sign of some coefficient in the normal form of the
fourth order is positive? And more generally: what is the orbits structure near the equilibrium
in this case, are there periodic orbits, invariant tori, their dimension, and so forth?

For the case when the equilibrium $O$ is of an elliptic-saddle-focus type, it has a
local smooth four-dimensional invariant center manifold $W^c$ corresponding to two pairs of
pure imaginary eigenvalues $\pm i(\omega_1\pm \sqrt{\eps_1})$ (one of two possible cases),
the restriction of the vector field on $W^c$ has an elliptic point at $O$, so there are two
one-parameter families of Lyapunov periodic orbits and we expect at some genericity assumptions
the existence of a positive measure set of
two-dimensional KAM tori in $W^c$. The center manifold is foliated into three-dimensional
submanifolds $W^c\cap \{H=const\}$, each KAM torus, if exists, belongs to some such
submanifold. The $W^c$ itself is normally hyperbolic in the transverse direction (because
four other eigenvalues have non-zeroth real parts). This implies the existence of the smooth
strong stable two-dimensional manifold $W^{ss}$ through $O$ and the smooth
strong unstable two-dimensional manifold $W^{uu}$. There are also two invariant six-dimensional
local center-stable $W^{cs}$ and center-unstable manifold $W^{cu}$ containing orbits which
approach to $W^c$ as $t\to \pm\infty$. Each 2-dimensional KAM torus $T\subset W^c$ belongs to
its own 7-dimensional level of the Hamiltonian $H=H(T)$ along with its 4-dimensional stable
$W^s(T)$ and unstable $W^u(T)$ manifolds. Hence again a question can be put: if these stable
and unstable manifolds intersect each other providing a (transverse in $H=H(T)$) homoclinic
orbit of $T$ with complicated nearby orbit dynamics or intersect stable (unstable) manifold of
neighboring tori leading to the phenomenon like the Arnold diffusion? To the questions of such
type one is rather hard to answer, see \cite{Kaloshin,GelfTur,Delshams}.

The case, when $O$ is an elliptic point, belongs to the KAM theory and we
expect here, under some non-degeneracy conditions, the existence of a
positive measure set of invariant four-dimensional KAM tori in a
neighborhood of $O$. Also, Lyapunov families of elliptic periodic orbits
can exist here as well as complicated mixture of saddle invariant tori of
intermediate dimensions with intersecting invariant manifolds and
invariant tori of lesser dimension than four being elliptic locally.

As the study of the Hamiltonian Hopf bifurcation shows \cite{Meer,GL}, for some
sign of a coefficient in the normal form of the fourth order the existence of homoclinic
orbits to $O$ is possible. So, we expect this be the case for the DHHB for some
combination of signs for coefficients $A,B,C$ in the normal form of the fourth order,
as well.

\section{Normal form of DHHB and its unfolding}

As we have seen, the quadratic normal form of the critical Hamiltonian
contains a semi-simple part, represented by two quadratic functions $S_1, S_2$,
$S_1=p_2q_1 - p_1q_2,~S_2=p_4q_3 - p_3q_4$, and the nilpotent part given by two functions
$N_1, N_2$, $N_1 = (q_1^2+q_2^2)/2,$ $N_2 = (q_3^2+q_4^2)/2,$ where
$\{S_1,S_2\}=0$, $\{N_1,N_2\}=0,$ $\{S_i,N_j\}=0,~i,j=1, 2$.
So, we should apply a normalization procedure for DHHB, similar to \cite{Meer}, where
the unperturbed (as $\eps_1 = \eps_2 = 0$) quadratic Hamiltonian is as follows
$$\begin{array}{l}
H_2^0= \omega_1(p_2q_1 - p_1q_2)+\frac12(q_1^2 + q_2^2)+
\omega_2(p_4q_3 - p_3q_4)+ \frac12(q_3^2 + q_4^2).
\end{array}
$$
We choose here both signs plus in front of $N_1, N_2$, that is not
essential for the study.

Now we introduce complex coordinates (here $i = \sqrt{-1}$)
\begin{equation}\label{comp}
\begin{array}{l}
q_1=\frac{1}{\sqrt{2}}(z_1-w_2),\;q_2=\frac{1}{i\sqrt{2}}(z_1+w_2),\;
p_1=\frac{1}{\sqrt{2}}(z_2+w_1),\;p_2=\frac{1}{i\sqrt{2}}(z_2-w_1),\;\\
q_3=\frac{1}{\sqrt{2}}(z_3-w_4),\;q_4=\frac{1}{i\sqrt{2}}(z_3+w_4),\;
p_3=\frac{1}{\sqrt{2}}(z_4+w_3),\;p_4=\frac{1}{i\sqrt{2}}(z_4-w_3),\;\\

dq_1\wedge dp_1+dq_2\wedge dp_2 +dq_3\wedge dp_3+dq_4\wedge dp_4 = \\
dz_1\wedge dw_1 + dz_2\wedge dw_2 + dz_3\wedge dw_3 + dz_4\wedge dw_4.
\end{array}
\end{equation}
In these symplectic complex coordinates we get
$$\begin{array}{l}
S_1 = i(z_1w_1+z_2w_2),\;S_2 = i(z_3w_3+z_4w_4),\;
N_1 = -z_1w_2,\;M_1 = \frac12(p_1^2+p_2^2) = z_2w_1,\\
N_2 = -z_3w_4,\;M_2 = \frac12(p_3^2+p_4^2) = z_4w_3,\\
H_2 = i\omega_1(z_1w_1+z_2w_2) - z_1w_2 + \eps_1z_2w_1 +
i\omega_2(z_3w_3+z_4w_4) - z_3w_4 + \eps_2z_4w_3.
\end{array}
$$

Let us first assume the absence of resonances between frequencies:
$\omega_1/\omega_2 \notin \mathbb Q$. Then we find the normal form with respect
to semi-simple part $S = \omega_1S_1+\omega_2S_2$, i.e. we should find those symplectic
transformations after which $\{H,S_1\}\equiv 0$ and $\{H,S_2\}\equiv 0$. This introduces
a symplectic $\mathbb{T}^2$-action generated by two commuting Hamiltonian flows $X_{S_1}$,
$X_{S_2}$. We notice, similar to \cite{Meer}, that the ring of functions being invariant
under this action is generated by the basic invariants
\begin{equation}\begin{array}{l}
S_1,\;S_2,\;N_1,\;N_2,\;
M_1 = \frac12(p_1^2+p_2^2),\;M_2 = \frac12(p_3^2+p_4^2),\\
P_1 = p_{1}q_{1}+p_{2}q_{2},\;P_2 = p_{3}q_{3}+p_{4}q_{4}\end{array}
\end{equation}
with the related Poisson structure relations given through Table 1. We compute the Poisson
brackets of any pair of invariants, called structure functions, and assemble them into a skew
symmetric matrix $G$, the structure matrix of the Poisson structure.
\begin{small}
\begin{center}
{
\renewcommand{\arraystretch}{1.3}
\begin{tabular}{c|cccccccc}
\multicolumn{8}{c}{Table 1: Poisson brackets for the invariants}\\ \hline
   $\{ \downarrow , \rightarrow \}$ & $S_1$ & $S_2$ & $N_1$ &
   $N_2$ & $M_1$ & $M_2$ & $P_1$ & $P_2$  \\ \hline

   $S_1$ & $0$ & $0$ & $0$ & $0$ & $0$ & $0$ & $0$ & $0$ \\

   $S_2$ & $0$ & $0$ & $0$ & $0$ & $0$ & $0$ & $ 0$& $0$  \\

   $N_1$ & $0$ & $0$ & $0$ & $0$ & $P_1$ & $0$ & $2N_1$ & $0$ \\

   $N_2$ & $0$ & $0$ & $0$ & $0$ & $0$ & $P_2$ & $0$  & $2N_2$\\

   $M_1$ & $0$ & $0$ & $-P_1$ & $0$ & $0$ & $0$ & $-2M_1$  & $0$ \\

   $M_2$ & $0$ & $0$ & $0$ & $-P_2$ & $0$ & $0$ & $0$  & $-2M_2$ \\

   $P_1$ & $0$ & $0$ & $-2N_1$ & $0$ & $2M_1$ & $0$ & $0$  & $0$ \\

   $P_2$ & $0$ & $0$ & $0$ & $-2N_2$ & $0$ & $2M_2$ & $0$  & $0$ \\
   \hline
\end{tabular}
}
\end{center}
\end{small}
All this shows that the normal form with respect to the semi-simple part
of the quadratic Hamiltonian is given as a series in invariant functions
(we assume the Hamiltonian be a real analytic or at least $C^\infty$-function in
some neighborhood of $O$)
\begin{equation}\label{semisimple}\begin{array}{l}
H = \omega_1S_1+\omega_2S_2 + N_1 + N_2  + \nu M_1 + \mu M_2 +\\
\sum_{i_1+i_2+j_1+j_2+k_1+k_2+l_1+l_2=2}^\infty
a_{i_1i_2j_1j_2k_1k_2l_1l_2}S_1^{i_1}S_2^{i_2}M_1^{j_1}M_2^{j_2}
N_1^{k_1}N_2^{k_2}P_1^{l_1}P_2^{l_2}.
\end{array}
\end{equation}
But we need to find the normal form w.r.t. the whole quadratic Hamiltonian $H_2 = H_2^s +
H_2^n$, not only w.r.t. the semi-simple part. So, we need to perform transformations for
which the normalized Hamiltonian $\mathcal H$ satisfies in addition to the identities
$\{\mathcal H,N_i\} \equiv 0,$ $i=1,2,$. Acting as in \cite{Meer} and taking into account
the nilpotent part, we come to the following theorem
\begin{theorem}
The normal form for the unfolding of the critical Hamiltonian \eqref{crit} expressed in terms
of the basic invariants reads as follows
\begin{equation}\label{NormalFormNilpotent}
\begin{array}{l}
\mathcal{H} = \omega_1S_1+\omega_2S_2 + N_1 + N_2  + \nu M_1 + \mu M_2 +
\sum_{i+j+k+l=2}^\infty A_{ijkl}S_1^iS_2^jM_1^kM_2^l,\;i,j,k,l\geq 0,
\end{array}
\end{equation}
where only invariants $S_i, M_i$, $i=1,2,$ are involved into the sum.
\end{theorem}

Explicitly, $\{S_1, \mathcal{H}\} = 0$, and $\{S_2, \mathcal{H}\} = 0$ and so $S_1$
and $S_2$ are integrals. We notice that in the space of invariants
$(S_1,S_2,N_1,M_1,P_1,N_2,M_2,P_2)$ two syzygies hold
$$4N_k M_k = P_k^2 + S_k^2 , N_k\geq 0, M_k \geq 0, k = 1,2\}.$$
Therefore, they have to be taken into account when studying orbits in this
space. Also, we conclude that the detuning parameters are embedded in $\nu,~\mu$.
As in the normal form for HHB \cite{Meer}, we limit ourselves, as a first
step, to the truncated normal form of the fourth order
\begin{equation}
\label{truncated-normal-form}
\begin{array}{l}
\bar{\mathbf{H}}^{\nu,\mu} = \omega_1S_1+\omega_2S_2 + N_1 + N_2 +\nu M_1 + \mu M_2 +
A M_1^2+ B M_1M_2+ C M_2^2,
\end{array}
\end{equation}
where $A=A_{0020},~B=A_{0011},~C=A_{0002}$ are renamed for simplicity.
\begin{remark}
It is important for the following study, the same normal form of the fourth order will
be for the unfolding, if we assume instead of irrationality $\omega_1/\omega_2$
the absence of strong resonances $\omega_1/\omega_2 \in \{1/4, 1/3, 1/2\}.$
\end{remark}

\subsection{The reduced phase space}

In any case, functions $S_1, S_2$ are integrals generating the Hamiltonian
action of the torus $T^2$ in $\R^8$. These two integrals are independent
at all points except points of two mutually complementary symplectic four-dimensional planes
$L_2 = \{q_1=q_2=p_1=p_2=0\}$ and $L_1 = \{q_3=q_4=p_3=p_4=0\}.$ Thus the common procedure
of reduction allows one to transform the system formally to the four-dimensional one
(two degrees of freedom). The related phase space of the reduced system will be either a smooth
four-dimensional symplectic manifold or a four-dimensional symplectic manifold with
singularities at some its points (an orbifold). Recall this procedure for our case, see details
in \cite{MW,CB}.

Consider a joint level of two integrals $S_1 = \zeta_1, S_2 = \zeta_2.$ This is the direct
product of related submanifolds in $L_1$ and $L_2,$ respectively. Near points where
differentials $dS_1, dS_2$  are independent (non-collinear for this case) the joint
level is a smooth 6-dimensional ball. For quadratic $S_i$, depending each
on only the mutually complementary half of coordinates, the dependence will occur only if one
of 4-dimensional co-vectors is zeroth, i.e. at the origin of the related 4-plane $L_i$. Thus,
the direct product of two 3-dimensional submanifolds situated within the
related 4-dimensional plane $L_i$ is as follows. Let us determine first these submanifolds.
For $\zeta_i = 0$ we get in $L_i$ the cone with the alone singularity at the origin,
topologically this is a cone over 2-torus\footnote{A cone over some topological
manifold $M$ is the topological space obtained from the space $I\times M$, $I=[0,1]$,
by the identification to the point of the subset $\{0\}\times M$.}. For $\zeta_i \ne 0$ one
gets a solid torus. To understand this, let us assume $\zeta_i > 0$ for definiteness, and
introduce coordinates $(x_1,x_2,y_1,y_2)$ in $L_i$, as follows:
$$
q_1=(x_1-y_2)/\sqrt{2},\;q_2=(y_1-x_2)/\sqrt{2},\;p_1=(x_2+y_1)/\sqrt{2},\;
p_2=(x_1+y_2)/\sqrt{2},\;
$$
The set is given by the equation $x_1^2 + x_2^2$ - $y_1^2 + y_2^2 = 2\zeta_i.$ If we go
around the circle $y_1^2 + y_2^2 = \rho^2$ we get in another plane the circle
$x_1^2 + x_2^2 = 2\zeta_i + \rho^2$. So, when $\rho = 0$ we get a circle in $L_i,$ but
if $\rho > 0$ we get a 2-torus in $L_i$. Varying $\rho$ in $[0,\infty)$ we have a solid
torus $S^1\times\R^2$. If $\zeta_i=0$, then at $\rho = 0$ we get a point
in $L_i$ but for $\rho > 0$ we have a 2-torus in $L_i.$ Varying $\rho$ from
zero to infinity we get a cone over a torus.

Now we should describe the products of two such 3-dimensional submanifolds
lying in the related 4-planes $L_1, L_2.$ Consider first the case of the
set given as $S_1 = \zeta_1$ and $S_2 = \zeta_2$ with $\zeta_1\zeta_2 \ne 0$.
This is the direct product of a solid torus in $L_1$ on a solid torus in
$L_2.$ If $\zeta_1 = 0$ or $\zeta_2 = 0$ but not both of them, then the set is the
direct of a cone over torus and a solid torus. This set contains singular
subset -- $\{0\}\times (S^1\times\R^2)$ -- at which points the set is not a
manifold, near any other point the set is a smooth 6-dimensional ball. If
both $\zeta_1 = \zeta_2 = 0$, then the set contains two singular subsets
glued at the common origin, near any other point the set is a smooth 6-dimensional
ball.

Through any nonsingular point of the set $S_1 = \zeta_1$ and $S_2 =
\zeta_2$ a two-dimensional torus passes being the orbit of the Hamiltonian
action of $T^2$. Through points of singular sets (except the origin) degenerate
one-dimensional orbits of the action pass (circles), the origin is the
singular point of the action. This implies that for nonzero
$\zeta_1\zeta_2$ the reduced phase space w.r.t. the action of $T^2$ we get
a smooth four-dimensional symplectic manifold, but when $\zeta_1 = 0$ or $\zeta_2 =
0$ we shall get symplectic orbifolds \cite{}.

\section{Truncated system via invariants}

The differential equations defined by the Hamiltonian (\ref{truncated-normal-form})
in invariants as variables are as follows:
$$
\dot M_i = \{M_i,H_4\},\;\dot N_i = \{N_i,H_4\},\;\dot P_i = \{P_i,H_4\},\;\dot S_i =
\{S_i,H_4\},\; i=1,2,
$$
or in the more explicit form, taking into account the bracket
relations of the Table 1, we come to the system
\begin{subequations}
\label{8eqs}
\begin{align}
\dot{M_1} & \;\; = \;\; -P_1, \\
\dot{M_2} & \;\; = \;\; -P_2, \\
\dot{N_1} & \;\; = \;\; P_1(2AM_1+BM_2+\nu), \\
\dot{N_2} & \;\; = \;\; P_2(2CM_2+BM_1+\mu), \\
\dot{P_1} & \;\; = \;\; 2M_1(2AM_1+BM_2+\nu)-2N_1, \\
\dot{P_2} & \;\; = \;\; 2M_2(2CM_2+BM_1+\mu)-2N_2, \\
\dot{S_1} & \;\; = \;\; 0,\\
\dot{S_2} & \;\; = \;\; 0.
\end{align}
\end{subequations}
The integrals $S_1=\zeta_1$ and $S_2=\zeta_2$ are considered as some
parameters, they enter to the system via the syzygies
\begin{equation}
\label{syzygies}
4N_i M_i - P_i^2 = \zeta_i^2,~~i=1,2,\;M_i, N_i\ge 0.
\end{equation}
We consider $8$-dimensional system (\ref{8eqs}) foliated into 6-dimensional subvarieties
defined by integrals and for each set of constants $\zeta_i,~~i=1,2,$ syzygies determine
invariant 4-dimensional submanifolds, we analyze the system on them. This
demonstrates the application of the reduction in invariants.

Thus, the reduction w.r.t. this action leads to a 2-DOF Hamiltonian system. When are these
reduced systems integrable or non-integrable, what are their structures, answering to these
questions are the primary goals for us.

Because $S_1, S_2$ are constants, the system reduces to the system
\begin{subequations}\label{6eqs}
\begin{align}
\dot{M_1} & \;\; = \;\; -P_1 \\
\dot{M_2} & \;\; = \;\; -P_2 \\
\dot{N_1} & \;\; = \;\; P_1(2AM_1+BM_2+\nu) \\
\dot{N_2} & \;\; = \;\; P_2(2CM_2+BM_1+\mu) \\
\dot{P_1} & \;\; = \;\; 2M_1(2AM_1+BM_2+\nu)-2N_1 \\
\dot{P_2} & \;\; = \;\; 2M_2(2CM_2+BM_1+\mu)-2N_2
\end{align}
\end{subequations}
with the phase space
$$
\{(M_1,M_2,N_1,N_2,P_1,P_2):~4N_i M_i = P_i^2 + \zeta_i^2,~~i=1,2,\;
M_i, N_i\ge 0,\;4N_i M_i \ge \zeta_i^2\}.
$$
This set can be considered as the product of the set in the positive ortant in $\R^4$
with coordinates $(M_1,M_2,N_1,N_2)$ defined by the inequalities $4N_i M_i \ge
\zeta_i^2$. For each point in this 4-dimensional set we have four points
defined by equalities $P_i = \pm\sqrt{4N_i M_i - \zeta_i^2},$ $i=1,2.$

In fact, first four equations of (\ref{6eqs}) in variables $(M_1, M_2, N_1, N_2)$ with
inserted $P_i = \pm \sqrt{4N_iM_i-\zeta_i^2},$ $i=1,2,$ define four-dimensional
systems depending on parameters $(\mu,\nu,\zeta_1,\zeta_2)$. Nevertheless, some calculations
are performed more conveniently in these coordinates, in particular, equilibria and their
types (see, Addendum).

The four-dimensional phase space can be described as follows.
Consider the direct product of two $\R^3_i, i=1,2,$ with coordinates
$(M_1,N_1,P_1)$ for $\R^3_1$ and $(M_2,N_2,P_2)$ for $\R^3_2.$ In each
$\R^3_i$ we choose the positive quadrant $M_i\ge 0, N_i\ge 0$ and in the
cylindric domain $D_i: 4M_iN_i - \zeta_i^2\ge 0$ we take a smooth (if $\zeta_i \ne 0$)
two-dimensional surface $C_i$ (a sheet of a hyperboloid) given as $P_i = \pm\sqrt{4M_iN_i -
\zeta_i^2}$. If $\zeta_i = 0,$ then this surface is a sheet of cone with a singular point at
the origin $(0,0,0)$ and this cone is smooth near any other point. Thus, for the fixed
constants $(\zeta_1, \zeta_2)$ in the six-dimensional space $\R^3_1\times\R^3_2$ we get
the direct product of the related two-dimensional manifolds (cones or hyperboloids) giving
a four-dimensional submanifold diffeomorphic to $\R^4$, if $\zeta_1\zeta_2 \ne 0$, or
a four-dimensional topological submanifold homeomorphic to $\R^4$, otherwise. This 4-manifold
is smooth at any point except for two-dimensional surface $\{(0,0,0)\}\times C_2$, if
$\zeta_1=0, \zeta_2 \ne 0,$ two-dimensional surface $C_1\}\times\{(0,0,0)\}$, if
$\zeta_2=0, \zeta_1 \ne 0,$ or two corresponding two-dimensional surfaces intersecting
only at the origin.

In fact, at any nonsingular point of the four-dimensional manifold defined by the
constants $\zeta_1, \zeta_2$ we get a coordinate chart defined by coordinates $(M_i,N_i,P_i),
i=1,2,$ and syzygies. The point on the manifold has four coordinates
$(M_1,N_1,M_2,N_2)$ along with the collection of four signs $(\pm,\pm)$, where the
first two signs select the sheet of the cone or hyperboloid for the first
quadrant $M_1\ge 0, N_1\ge 0$ and the second two signs related with sheets
in the second quadrant $M_2\ge 0, N_2\ge 0.$ If the projection of a smoothness point belongs
to the branching line of the first quadrant or the second quadrant or both of
them, then we take near this point as coordinates either $(M_1,P_1,M_2,N_2)$ or
$(N_1,P_1,M_2,N_2)$ and the sign corresponding to the second pair. Similarly we do for
the other cases.

The system (\ref{8}) provides uniquely the related vector field near this point in coordinates.
A transition from one coordinate system to another one is determined by
the movement along the related orbit. The transition can arise, if the
projection of the orbit crosses the branching line of the related quadrant.
If it occurs at the nonsingular point of the branching line, then the projection of the orbit
to the related quadrant makes the passage onto the another sheet with changing sign of
related $P_i.$

It remains to understand how an orbit continues when some of its
projections on the related quadrant tends to the origin. This means, as a corollary,
that the related $\zeta_i$ is zero and the projection lies on the related
cone and we approach in the 4-manifold to the invariant two-dimensional singular
submanifold. The behavior of its orbits are described in Subsection 4.1, so that the orbit will
move near some orbit of this submanifold.

At any smoothness point of the 4-dimensional manifold obtained the
system (\ref{6eqs}) generates a Hamiltonian vector field. Actually, we are
interested in the orbit structure of this vector field in some
neighborhood of its equilibrium at the origin.

First notice that due to our assumption $4AC>B^2$ the signs of $A$ and $C$ are same.
The equilibria are solutions of the system:
$$\begin{array}{l}
P_1 = 0,\;P_2=0,\;N_1= M_1(2AM_1+BM_2+\nu),\;N_2=M_2(BM_1+2CM_2+\nu),\\
M_1^2(2AM_1+BM_2+\nu)=0,\; M_2^2(BM_1+2CM_2+\mu)=0.
\end{array}
$$
From two last equations w.r.t. $M_1, M_2$ we get the following collection of
equilibria:\\
$(M_1,M_2,N_1,N_2,P_1,P_2) = $
\begin{enumerate}
\item $E_1 = (0,0,0,0,0,0)$;
\item $E_2 = (0,-\frac{\mu}{2C},0,0,0,0)$;
\item $E_3 = (-\frac{\nu}{2A},0,0,0,0,0)$;
\item $E_4 = (\frac{B\mu-2C\nu}{4AC-B^2},\frac{\nu B - 2A\mu}{4AC-B^2},0,0,0,0).$
\end{enumerate}
Because parameters $\mu,\nu$ can take any signs, the equilibria $E_2-E_4$ exist not
always, but only when signs of $M_1, M_2$ are nonnegative. Their types
will be described at the Addendum.

\subsection{Systems on invariant 3-planes}

The system (\ref{6eqs}) has two invariant 3-planes: $M_1 = N_1 = P_1 = 0$,
$M_2 = N_2 = P_2 = 0$. On the first 3-plane one gets $\zeta_1=0$, on the second 3-plane
one gets $\zeta_2=0$. Each of 3-planes is foliated into invariant subvarieties
$4M_2N_2 = P_2^2 + \zeta_2^2$ (the first one) and $4M_1N_1= P_1^2 + \zeta_1^2$ (the second one).
On the second plane with coordinates $(M_1,N_1,P_1)$ we come, with the account of the related
syzygy, to two-dimensional systems depending on parameters $\nu,\zeta_1$
\begin{equation}\label{3eq}
\begin{array}{l}
\dot{M_1} = -P_1,\;4M_1N_1 = P_1^2 + \zeta_1^2,\\
\dot{N_1} = P_1(2AM_1+\nu),\\
\dot{P_1} = 2M_1(2AM_1+\nu)-2N_1,
\end{array}
\end{equation}
and similarly for the second 3-plane with the change $(M_1,N_1,P_1)$ to
$(M_2,N_2,P_2)$ and parameters $\nu,A,\zeta_1$ to $\mu,C,\zeta_2$. The system (\ref{3eq})
has two integrals $\nu M_1 + N_1 + AM_1^2$ and $4M_1N_1 - P_1^2$, as it follows from the
first two equations, so its phase space is foliated into levels $\nu M_1 + N_1 + AM_1^2 = k$
being parabolic cylinders in 3-space parallel to the $P_1$-axis\footnote{This is a
ruled surface with the directrix being the parabola $N_1 = k-\nu M_1 - AM_1^2$ and all its
straight-lines parallel to $P_1$-axis.}. So, the orbits of the system with
fixed parameters $\nu, A$ are obtained by the intersection of the surface $4M_1N_1 =
P_1^2 + \zeta_1^2$ at fixed $\zeta_1$ with cylinders $\nu M_1 + N_1 + AM_1^2 = k$ as $k$
varies. The related foliation of the cone into curves depends on the sign
of $A$, the value of $\nu$ is a small bifurcation parameter.
This study is similar to what was done in \cite{Hanssmann}, Sec.2.2 and \cite{KuLe},
Sec.3.3.

\begin{figure}[htp] 
\centering
\subfigure[$A > 0,\nu > 0$]{%
\includegraphics[width=0.4\textwidth]{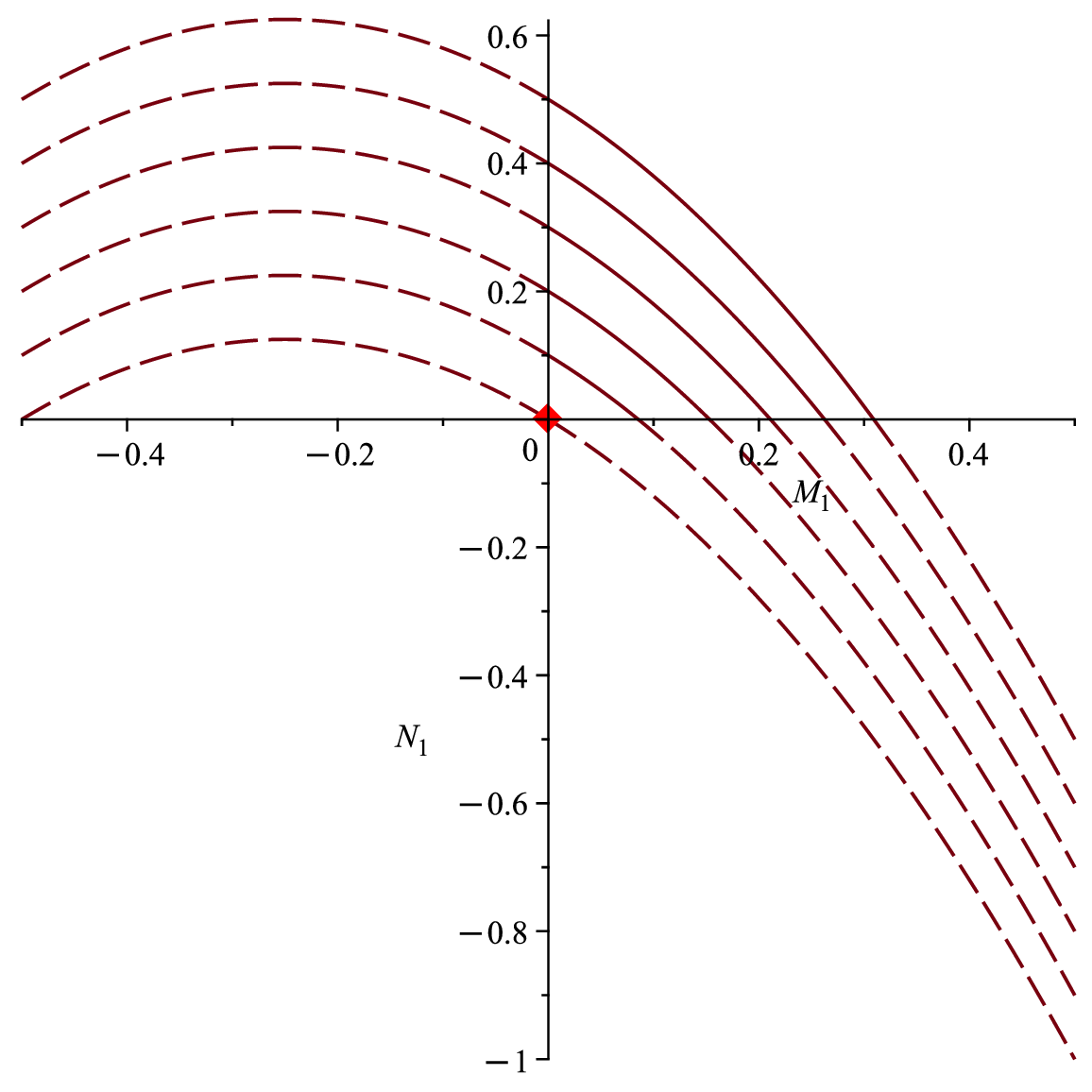}%
\label{fig:1:a}%
}\hfil
\subfigure[$A > 0,\nu < 0$]{%
\includegraphics[width=0.4\textwidth]{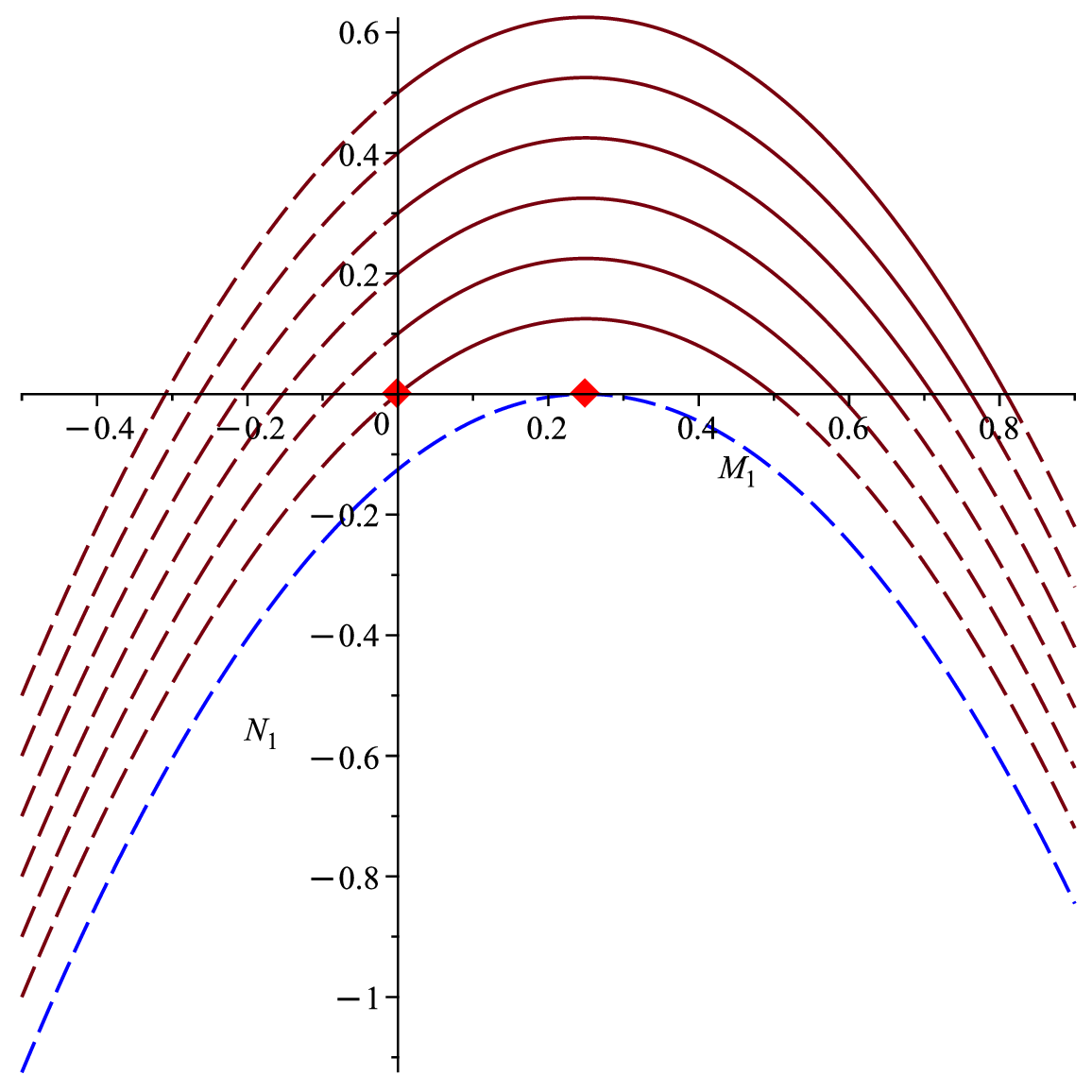}%
\label{fig:1:b}%
}

\subfigure[$A < 0,\nu > 0$]{%
\includegraphics[width=0.4\textwidth]{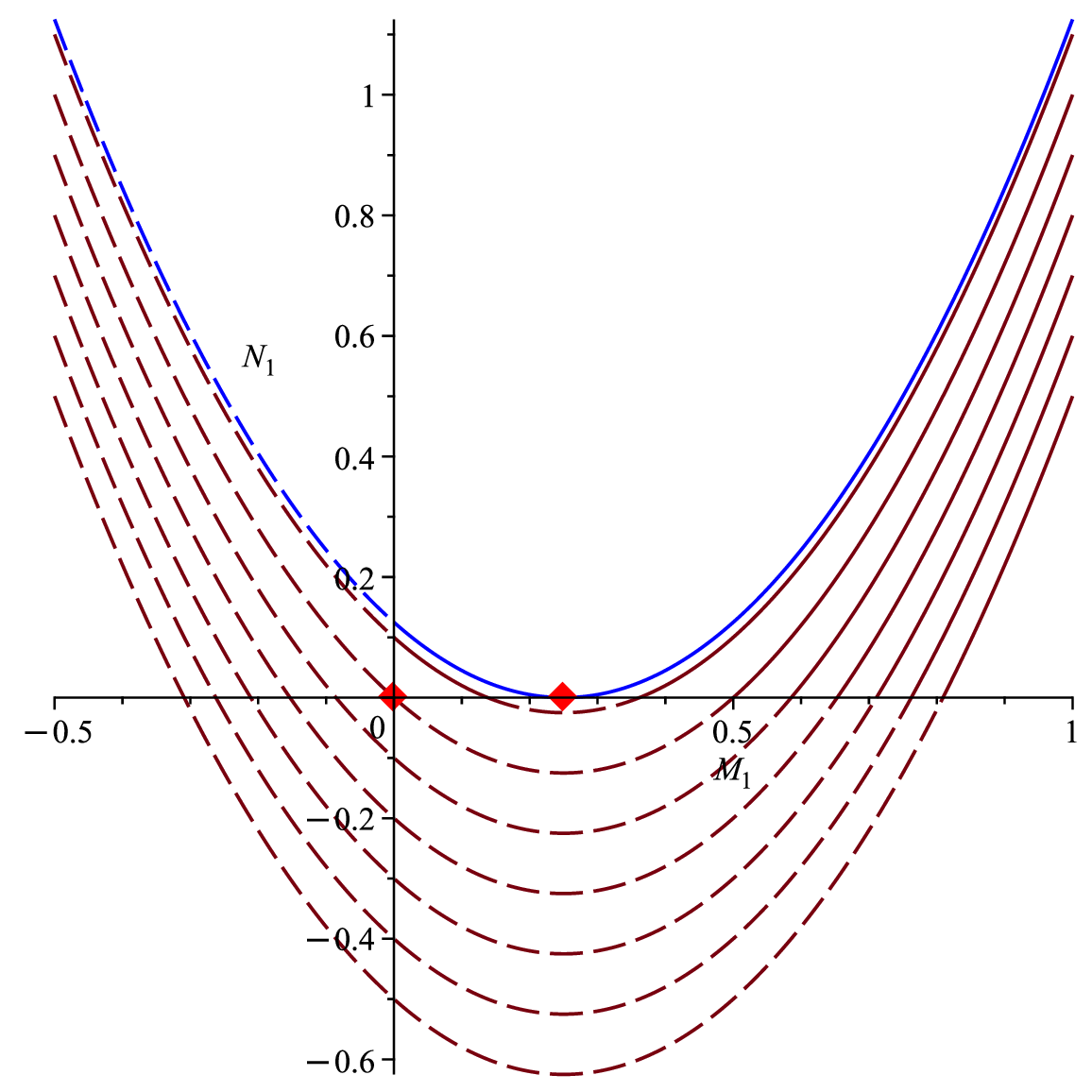}%
\label{fig:1:c}%
}\hfil
\subfigure[$A < 0,\nu < 0$]{%
\includegraphics[width=0.4\textwidth]{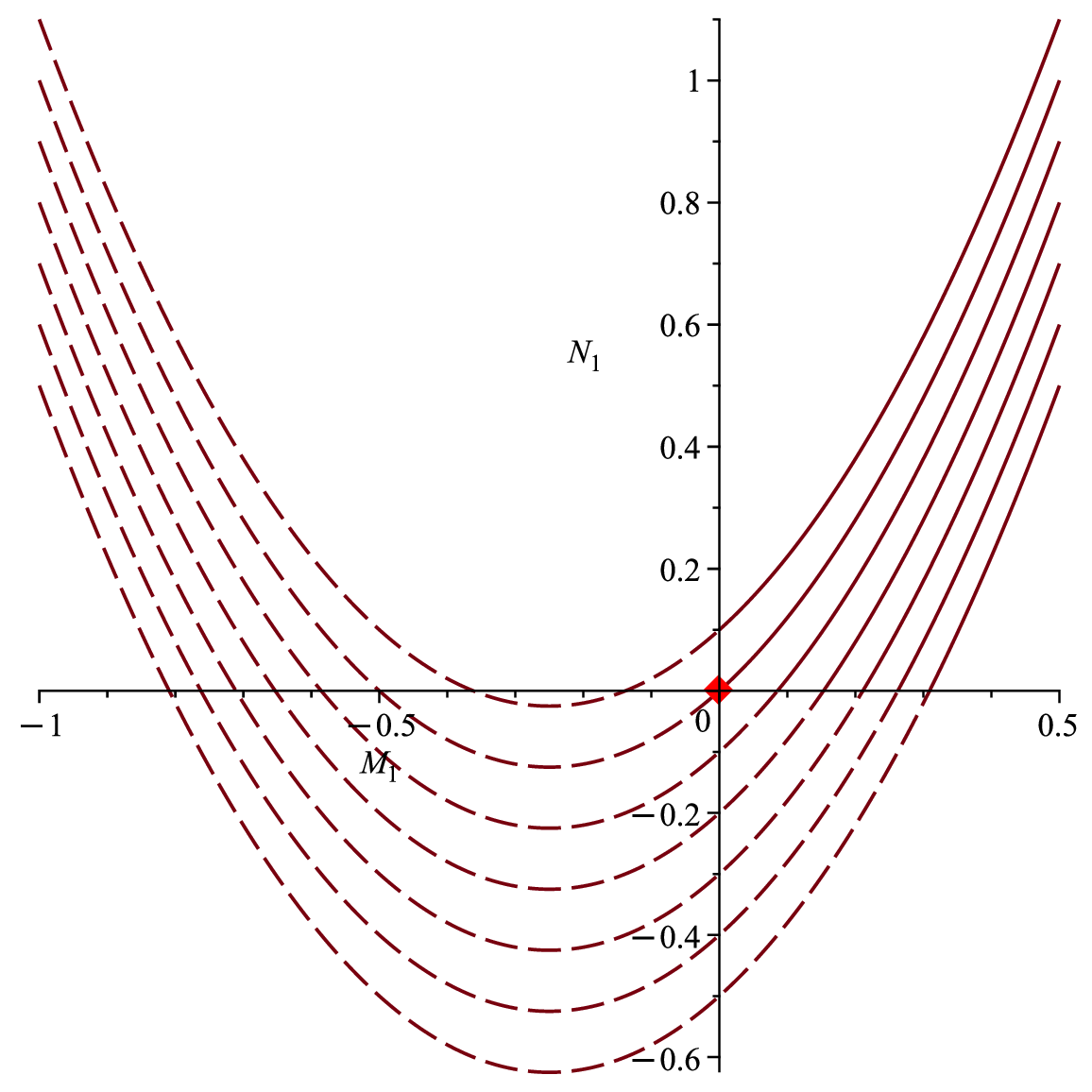}%
\label{fig:1:d}%
}
\caption{The of orbits on the cone in different cases: Indeed in these cases, we choose
the parameters $A =\pm 2$, and $\nu = \pm 1$ and plot the projection of orbits on
$(M_1,N_1)$-plane, where the projection of orbits on the cone are shown by solid lines,
while those out of the cone are dashed lines. The red points are projections of equilibria.
The blue line distinguishes by including tangency points in phase space.}
\label{fig:1}
\end{figure}

{\bf The case $\zeta_1=0$}. In the 3-space $(M_1,N_1,P_1)$ we obtain a smooth (except for the point
$(0,0,0)$) two-sheeted surface over the first quadrant $M_1\ge 0, N_1\ge 0$ given by
the relation $P_1 = \pm 2\sqrt{M_1N_1}$ (this is a sheet of the cone). The
branching line (the line of folds w.r.t. the projection onto $(M_1,N_1)$-plane)
is the boundary of the first quadrant $M_1\ge 0, N_1=0$ and
$N_1\ge 0, M_1=0.$ The tangent planes to the cone at these points (except
the non-smoothness point $(0,0,0)$) are vertical planes $N_1=0$ and $M_1=0,$
respectively. The orbit of (\ref{3eq}) through a point on the cone is the intersection
of the cone with that parabolic cylinder $\nu M_1 + N_1 + AM_1^2 = k$ which passes through
this point. At fixed $\nu, A$ and varied $k$ all these parabolic cylinders are parallel
to each other with respect to the shift along $N_1$-axis. Some orbits can reach the branching
line from upper $P_1>0$ or lower $P_1<0$ half of the cone, then they continue moving to
opposite part of the cone. The tangent to the orbit at the branching point is vertical
(parallel to $P_1$-axis). So, in order to get an image of the foliation of the cone
into orbits of (\ref{3eq}) one needs to take two copies of the related foliation on
the quadrant $M_1\ge 0, N_1\ge 0$ (see Fig.~\ref{fig:1}) and glue them along
the boundary.

There is the only parabolic cylinder of the $k$-family which is tangent to the cone, it
contains the point $(-\nu/2A,0,0)$ on the $M_1$-axis, i.e. at $k = -\nu^2/4A$. The
corresponding

point on the invariant cone is an equilibrium of the system (\ref{3eq}), it exists only
if $-\nu/2A > 0$. This equilibrium is a center, if $A>0$ and a saddle if $A<0$. Indeed, the cone
near this point is represented as a graph of the function $N_1 = P_1^2/4M_1$. Inserting
this into the third equation in (\ref{3eq}) we get the system of two differential equations
in coordinates $(M_1,P_1)$ (see (\ref{hyp}) below for the general case) with the equilibrium
at the point $(M_1,P_1) = (-\nu/2A,0)$ and its characteristic equation
(\ref{char}). For $\nu < 0, A>0$ the term $4(3AM_1+\nu) = -2\nu > 0$, so we have the center,
parabolic cylinders with $-\nu^2/4A < k < 0$ intersect the cone along closed orbits enclosing
the center. As $k\to -0$, the periodic orbit approaches to the closed curve being the separatrix
loop of the equilibrium at the vertex of the cone $(0,0,0).$ When $k$ increases further
becoming positive, the related parabolic cylinders give again periodic orbits on
the cone.

For $A>0,\nu > 0$ parabolas on the plane $(M_1,N_1)$, being directrices of parabolic
cylinders from $k$-family, have their vertices on the line $M_1 = -\nu/2A \le 0$, hence
for $k>0$ they intersect of a parabolic cylinder with the cone give each a unique periodic
orbit that contracts to the equilibrium $(0,0,0)$ as $k\to +0$, this equilibrium is of
the center type (but this is the non-smoothness point of the cone!). See Fig. \ref{fig:1:a}.

How to see the gluing of two pieces of the parabolas lying on the upper sheet of the surface
and on its lower sheet at their boundary points forming a smooth curve. Indeed, the boundary
point on the $N_1$-axis has coordinates $(0,N^b_1,0),$ $N^b_1>0,$ and the corresponding vector
of the system (\ref{3eq}) is vertical: $(0,0,-2N_1^b)$. So, upper and lower
parts of the orbits on two sheets with the same projection being the
pieces of parabola are glued smoothly.

Let now $A$ be negative. Since we have again the quadrant $M_1\ge 0,$ $N_1\ge 0,$
then the related parabolas $\nu M_1 + N_1 + AM_1^2 = k$ will have their wings going
to $\infty$ as $|M_1|\to \infty.$ Thus, for $\nu \le 0$ the line of parabola
vertices on the plane $(M_1,N_1)$ projects at the point $(-\nu/2A,0$ and the intersection
of of $k$-family of parabolic cylinders with the cone gives an infinite orbit on the cone
going to infinity in both directions of time  (see Fig.~\ref{fig:1}, downer row). For $k=0$ this curve passes through the vertex
$(0,0,0)$, hence this gives two separatrices of the equilibrium at the point $(0,0,0).$

If $\nu > 0$, the situation is different. In this case the equilibrium
at $(-\nu/2A,0,0)$ on the cone is a saddle, since free term of the characteristic equation
(\ref{char}) is negative. The parabolic cylinder through this point has
two components of the intersection with the cone. One component has
compact closure and form a separatrix loop of the saddle. Orbits on the
cone lying inside of this loop are closed and shrink as $k\to 0$ to the center at the
vertex of the cone. Another component is infinite and form two other separatrices of
the saddle going to infinity in both direction of time.

\begin{figure}[htp] 
\centering
\subfigure[$A > 0,\nu > 0$]{%
\includegraphics[width=0.4\textwidth]{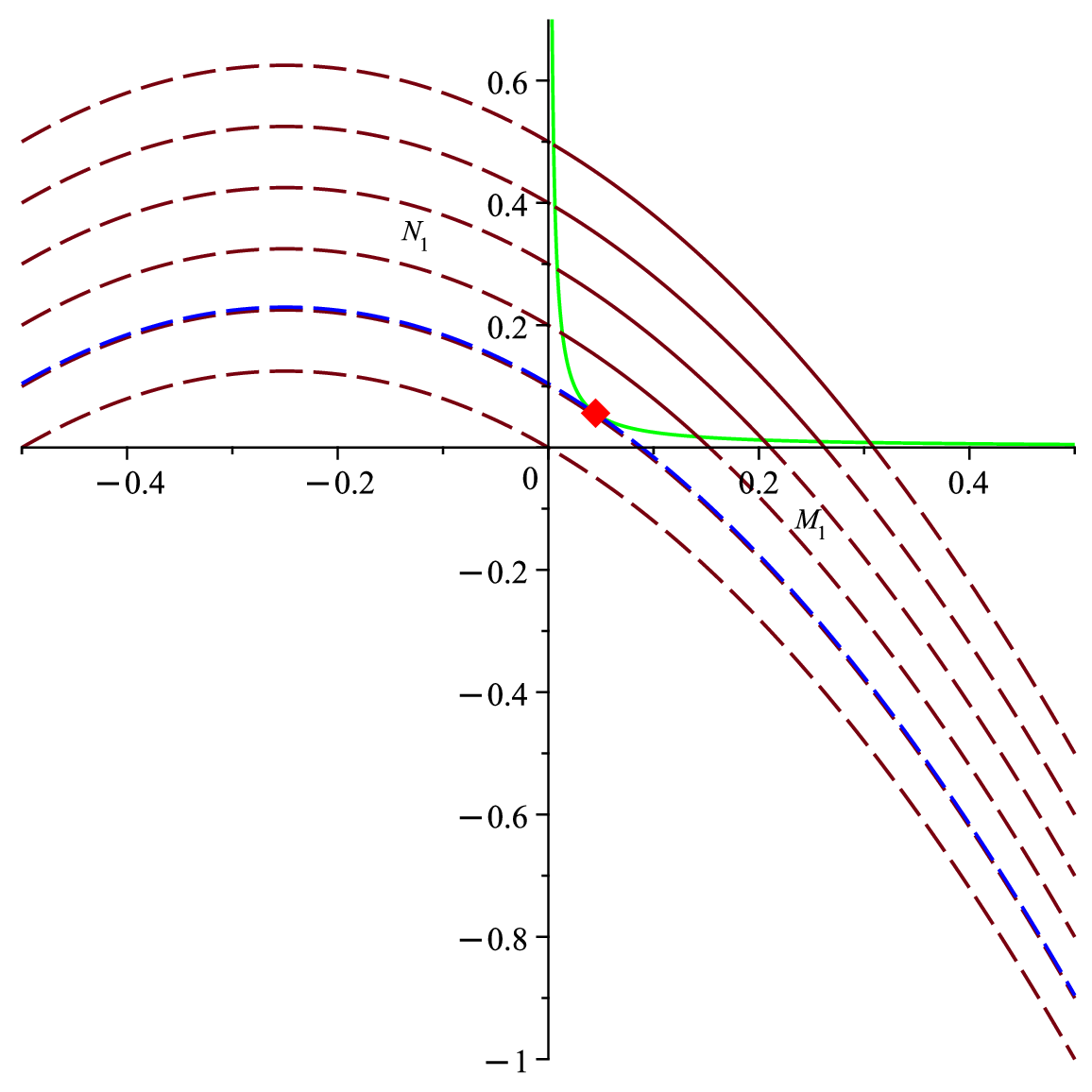}%
\label{fig:2:a}%
}\hfil
\subfigure[$A > 0,\nu < 0$]{%
\includegraphics[width=0.4\textwidth]{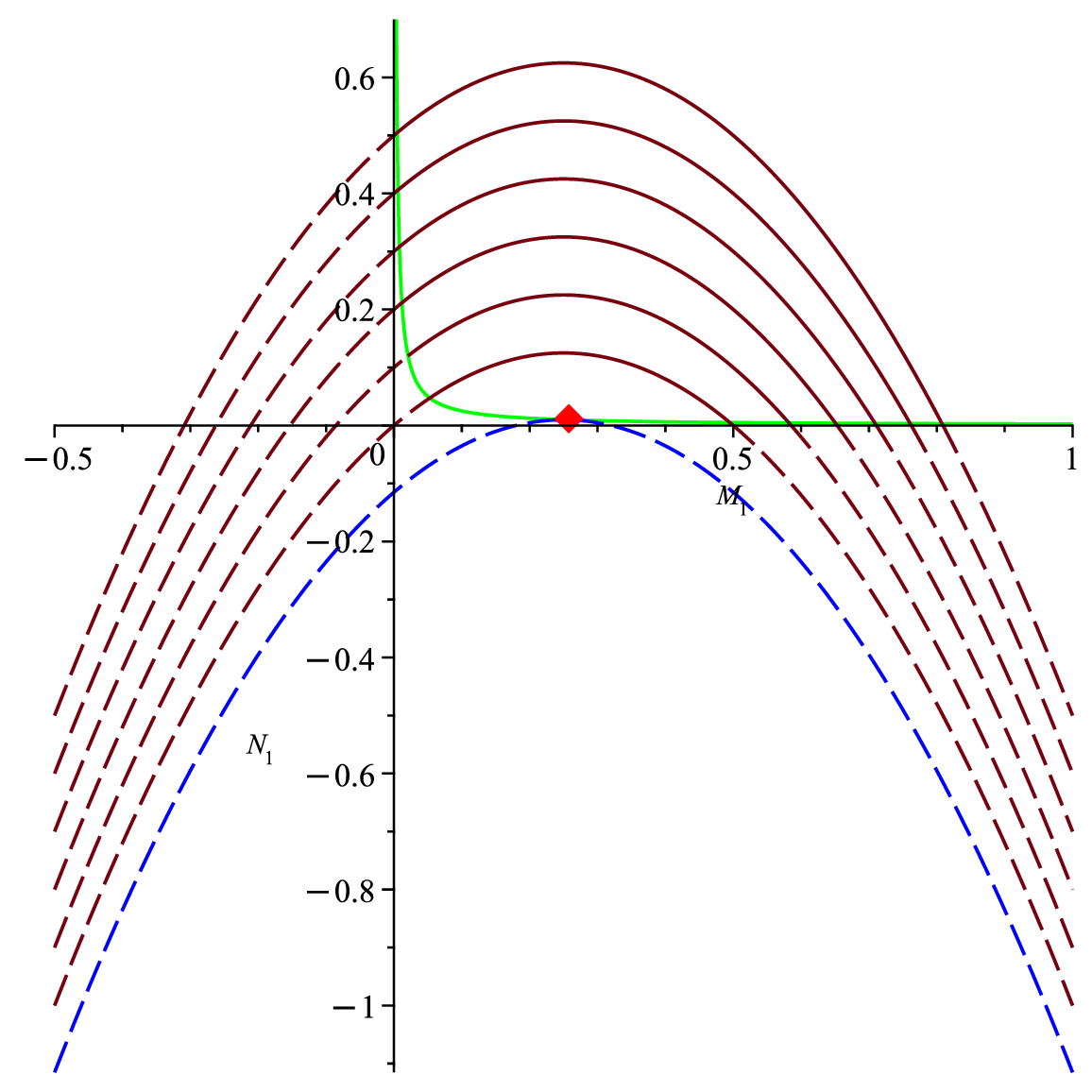}%
\label{fig:2:b}%
}

\subfigure[$A < 0,\nu > 0$]{%
\includegraphics[width=0.4\textwidth]{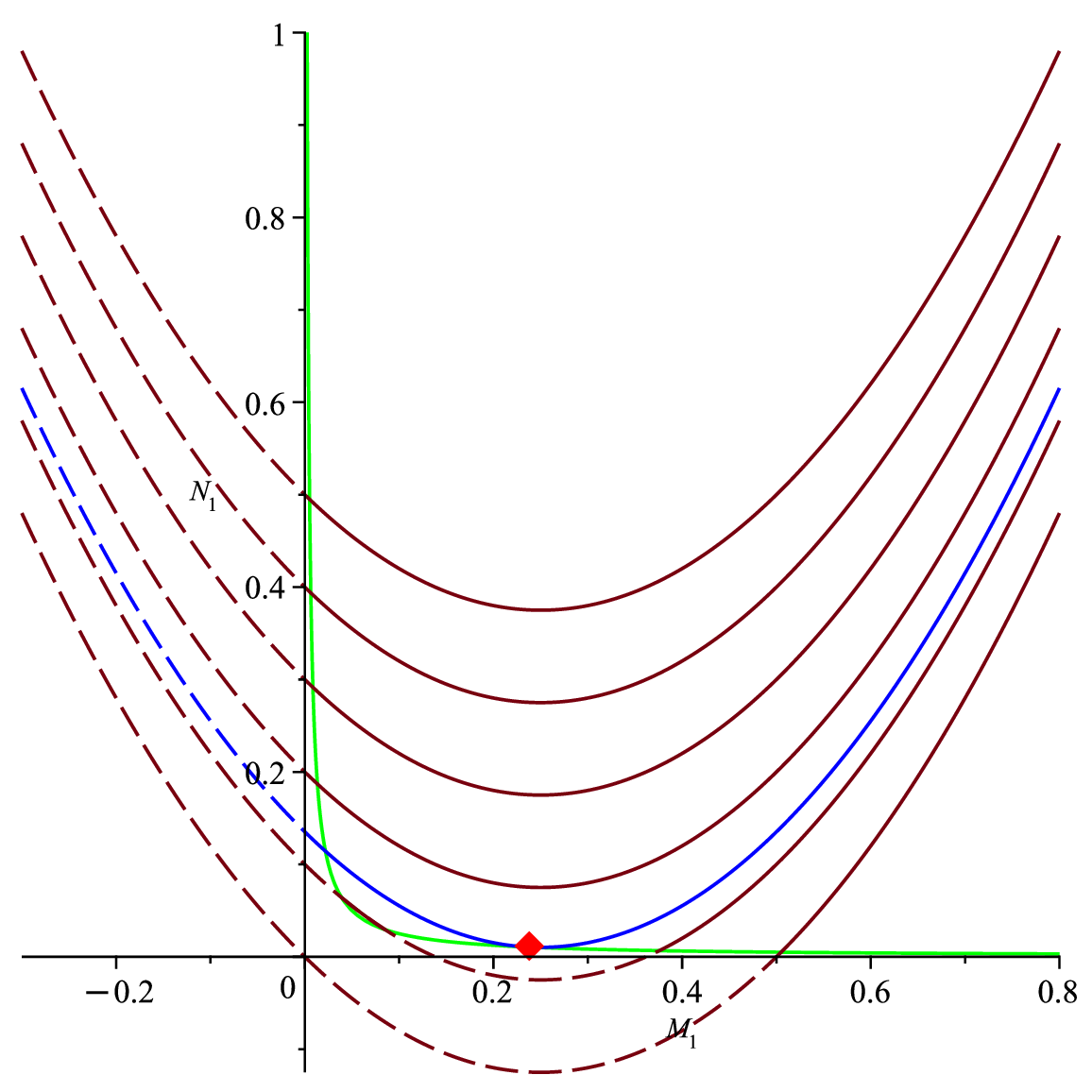}%
\label{fig:2:c}%
}\hfil
\subfigure[$A < 0,\nu < 0$]{%
\includegraphics[width=0.4\textwidth]{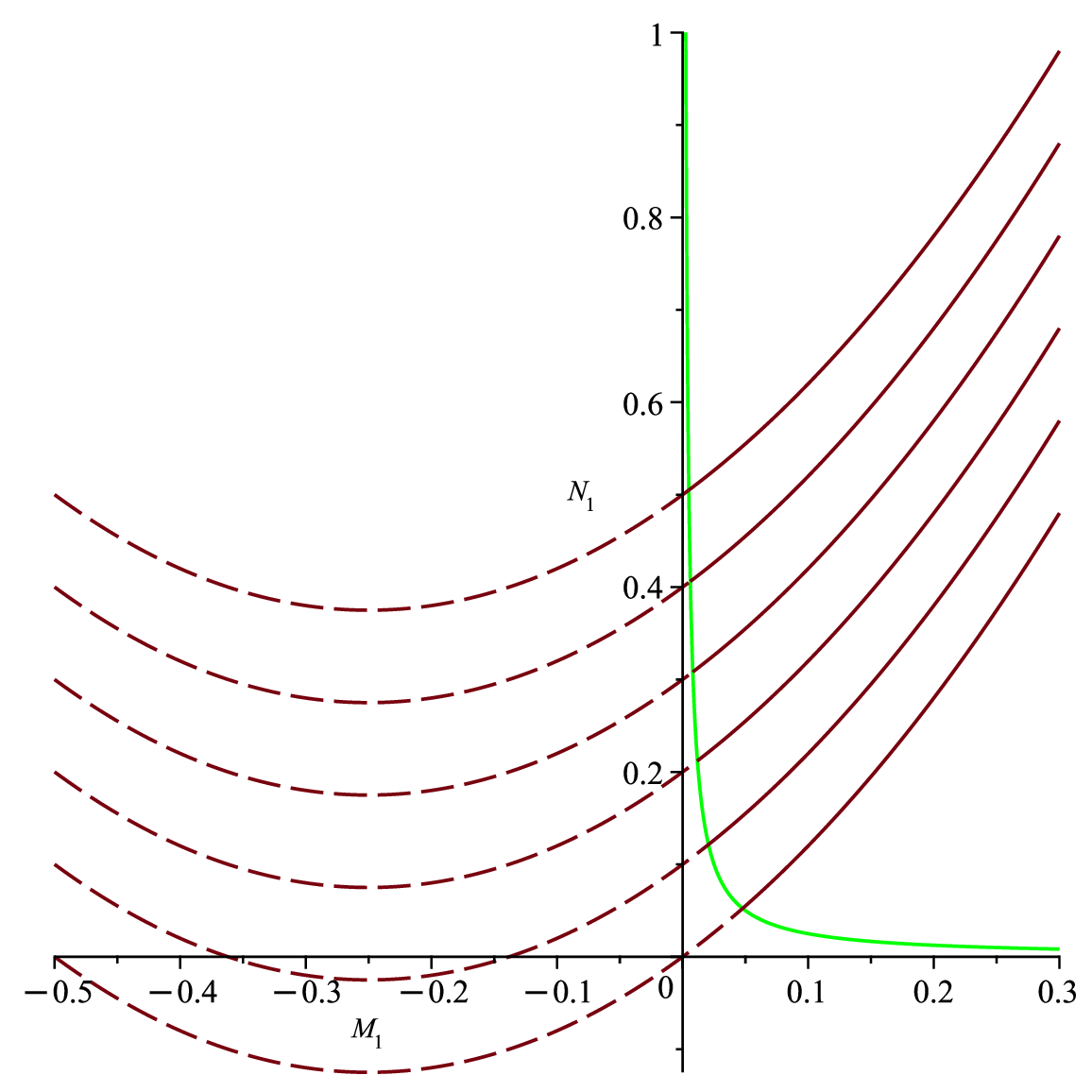}%
\label{fig:2:d}%
}
\caption{The orbits on the hyperboloid in different cases: Indeed in these cases,
we choose the parameters $A =\pm 2$, $\nu = \pm 1$ and $\zeta_1=\frac{1}{10}$ and plot
the projection of orbits on $(M_1,N_1)$-plane, where the projection of orbits on the hyperboloid
are shown by solid lines, while those out of hyperboloid are dashed lines. The red points are
projections of equilibria. The blue line distinguishes by including tangency points in
phase space.}
\label{fig:2}
\end{figure}

{\bf The case $\zeta_1 \ne 0$}. Here we have a surface being smooth at all
its points, it lies over a domain in the first quadrant separated out through the inequality
$4M_1N_1 \ge \zeta_1^2.$ This surface is the one sheet of two-sheeted hyperboloid
$4M_1N_1-P_1^2 = \zeta_1^2$ lying over the first quadrant. Over any point $(M_1, N_1)$ of
the open part of the domain we have two points of the hyperboloid defined as $P_1 = \pm \sqrt{4M_1N_1 -
\zeta_1^2}$, these two points coincide on the boundary -- branching line $4M_1N_1 = \zeta_1^2$
being the line of folds of the hyperboloid w.r.t. the projection on the $(M_1,N_1)$-plane.
The tangent planes at the points of folds are vertical planes. Here again in order to get
an image of the foliation of the hyperboloid into orbits of (\ref{3eq}) one needs to take
two copies of the related foliation on the quadrant $M_1\ge 0, N_1\ge 0$ (see Fig.~\ref{fig:2})
and glue them along the boundary.

As in the case $\zeta_1=0,$ we have the invariant parabolic cylinders $N_1 = k - \nu M_1 - AM_1^2$,
but now we consider their intersections with the hyperboloid given by a fixed $\zeta_1\ne 0$.
The equilibria of the system (\ref{3eq}) on the hyperboloid are its tangency points with
the parabolic

cylinders, so these equilibria belong to the branching line and are determined by
the system
$$
P_1=0, \;8AM_1^3 + 4\nu M_1^2 - \zeta_1^2 = 0,\;N_1 = \zeta_1^2/4M_1.
$$
The number of roots of the cubic equation depends on constants $\nu,\zeta_1$ and the sign
of $A.$ The equation can have either one or three simple roots in dependence of a location
of a point $(\nu,\zeta_1)$ in this parameter plane w.r.t the semi-cubic parabola
$4\nu^3 - 27A^2\zeta^2 = 0,$ but only roots with $M_1\ge 0$ are relevant.
As we shall see, the cases $A>0$ and $A<0$ are different.

Let $A$ be positive. Then the parabolic cylinders can be tangent with the hyperboloid at
one point with positive $M_1$, since for $A>0$ the equation $f(M_1)= 8AM_1^3 + 4\nu M_1^2 -
\zeta_1^2 = 0$
has the only positive root. Indeed, the derivative $8M_1(3AM_1+\nu)$ has zeroth at $M_1 = 0$
and $M_1^* = - \nu/3A$, if $\nu \ne 0$. At negative $\nu$ the point $M_1^* > 0$ is the minimum
of the cubic polynomial and $M_1=0$ is its maximum, so we get $f(M_1^*) < f(0) = -\zeta^2 <0.$
For positive $\nu$ one has the opposite case: the point $M_1^* < 0$ is its maximum and $M_1=0$
is the minimum, so $f' > 0$ for all $M_1\ge 0$. Thus, we conclude that for $A>0$ and any $\nu$
there is the only tangency points of the hyperboloid and parabolic
cylinders (see Fig.~\ref{fig:2}, upper row).

Passing again to the representation of the hyperboloid near the equilibrium to coordinates
$(M_1,P_1)$ we get the system of two differential equations
\begin{equation}\label{hyp}
\dot M_1 = -P_1,\;\dot P_1 = -\frac{P_1^2+\zeta_1^2}{2M_1} + 4AM_1^2 + 2\nu M_1 =
(8AM_1^3 + 4\nu M_1^2 -\zeta_1^2 - P_1^2)/2M_1.
\end{equation}
Let us denote $f(M_1) = 8AM_1^3 + 4\nu M_1^2 -\zeta_1^2$. For the linearization at
the equilibrium, where $f(M_1)=0, P_1=0$, the characteristic equation is as
follows
\begin{equation}\label{char}
\lambda^2 + (M_1f'-f)/2M_1^2 = \lambda^2 + 4(3AM_1+\nu) = 0,
\end{equation}
which has for positive root $M_1$, $A>0, \nu>0$ two imaginary eigenvalues, because for
the root $M_1$ we have $4(3AM_1+\nu) > 0$. Thus, if $A > 0$ the unique positive root
corresponds to the center type equilibrium on the hyperboloid and all other orbits are closed.

Let $A$ be negative. Then the same cubic equation for $\nu > 0$ can have, when
 $4\nu^3 > 27A^2\zeta^2$, two positive roots, separated by the critical point
$M_1 = -\nu/3A > 0,$ and no positive roots, if $0 < 4\nu^3 < 27A^2 \zeta^2$. For $\nu < 0$
the equation has not positive roots, since the derivative $8M_1(3AM_1+\nu)$ is negative
for positive $M_1$. Doing the same as for $A>0$, we come to the same differential system
in variables $(M_1,P_1)$ near the equilibria with the same characteristic equation.
The equilibrium corresponding to left positive root $M_1 + \nu/3A < 0,$ is a center and
that to the right $M_1 + \nu/3A > 0$ is a saddle. Related foliation when projecting
on the plane $(M_1,N_1)$ are on the Fig.~\ref{fig:2}.

\section{Addendum: other equilibria}

\subsection{The case of a saddle-saddle at the origin}

Let $\nu<0,~\mu<0$, then the equilibria of the system \eqref{6eqs} are as follows
\begin{itemize}
\item $E_1=(M_1,M_2,N_1,N_2,P_1,P_2)=(0,0,0,0,0,0)$ with eigenvalues $\pm2\sqrt{-\mu},
~\pm2\sqrt{-\nu},0,0$, that shows the origin is of type saddle-saddle-fold. So on a invariant
codimension 2 (4-dimensional) manifold origin is of type saddle-saddle for reduced system.
\item If $C>0$, $E_2=(N_1,N_2,M_1,M_2,P_1,P_2)=(0,0,0,\sqrt{-\frac{\mu}{2C}},0,0)$ then if
$B\mu - 2C\nu>0$, with eigenvalues $\pm \sqrt{-2\mu}i,~\pm \sqrt{\frac{2(B\mu - 2C\nu)}{C}},0,0$,
that shows it is of type saddle-center-fold and if $B\mu - 2C\nu<0$, with eigenvalues
$0,0,\pm \sqrt{-2\mu}i,~\pm \sqrt{-\frac{2(B\mu - 2C\nu)}{C}}i$, that shows it is of type
fold-center-center.
\item If $A>0$, $E_3=(N_1,N_2,M_1,M_2,P_1,P_2)=(0,0,\sqrt{-\frac{\nu}{2A}},0,0,0)$,
then if $B\nu - 2A\mu>0$, with eigenvalues $\pm \sqrt{-2\nu}i,~\pm \sqrt{\frac{2(B\nu - 2A\mu)}
{A}},0,0$, that shows it is of type saddle-center-fold and if $B\nu - 2A\mu<0$, with
eigenvalues $\pm \sqrt{-2\nu}i,~\pm i\sqrt{-\frac{2(B\nu - 2A\mu)}{A}},0,0$, that shows
it is of type center-center-fold.
\item If $B\nu - 2A\mu>0$, and $B\mu - 2C\nu>0$, then $$E_4=(N_1,N_2,M_1,M_2,P_1,P_2)
= (0,0,\sqrt{\frac{B\mu - 2C\nu}{4AC-B^2}},\sqrt{\frac{B\nu - 2A\mu}{4AC-B^2}},0,0)$$ where
the eigenvalues are the roots of quartic polynomial $\alpha x^4 + \beta x^2 +\gamma$, with
\begin{eqnarray*}
\alpha &=& 4AC-B^2>0,\\
\beta &=& 4A(B\mu - 2C\nu)+4C(B\nu - 2A\mu),\\
\gamma &=&  4A(B\mu - 2C\nu)(B\nu - 2A\mu)>0,\\
\end{eqnarray*}
and so if $A>0$ and $C>0$, then $E_4$ is of type hyperbolic-hyperbolic-fold, and if $A<0$
and $C<0$ with $A \,B^{2}+4 A B C +4 A \,C^{2}-2 B^{3}<0$ then $E_4$ is of type
saddle-saddle-fold or center-center-fold, and if $A<0$ and $C<0$ with $A \,
B^{2}+4 A B C +4 A \,C^{2}-2 B^{3}>0$ then $E_4$ is of type hyperbolic-hyperbolic-fold.
\end{itemize}

\subsection{The saddle-center case at the origin}

This case deals with either $\nu>0,~\mu<0$, or $\nu<0,~\mu>0$. So if $\nu>0,~\mu<0$, then
the equilibria of \eqref{6eqs} are as follows,
\begin{itemize}
\item $E_1=(N_1,N_2,M_1,M_2,P_1,P_2)=(0,0,0,0,0,0)$ with eigenvalues $\pm2\sqrt{-\mu},~\pm2\sqrt{\nu}i,0,0$,
that shows the origin is of type saddle-center-fold.
\item If $C>0$, $E_2=(N_1,N_2,M_1,M_2,P_1,P_2)=(0,0,0,\sqrt{-\frac{\mu}{2C}},0,0)$ then
if $B\mu - 2C\nu>0$, with eigenvalues $\pm \sqrt{-2\mu}i,~\pm \sqrt{\frac{2(B\mu - 2C\nu)}{C}},0,0$,
that shows it is of type saddle-center-fold and if $B\mu - 2C\nu<0$, with eigenvalues
$\pm \sqrt{-2\mu}i,~\pm \sqrt{-\frac{2(B\mu - 2C\nu)}{C}}i,0,0$, that shows it is of type
center-center-fold.
\item If $A<0$, $E_3=(N_1,N_2,M_1,M_2,P_1,P_2)=(0,0,\sqrt{-\frac{\nu}{2A}},0,0,0)$ then if
$B\nu - 2A\mu>0$, with eigenvalues $\pm \sqrt{2\nu},~\pm \sqrt{-\frac{2(B\nu - 2A\mu)}{A}}i,0,0$,
that shows it is of type saddle-center-fold and if $B\nu - 2A\mu<0$, with eigenvalues
$\pm \sqrt{2\nu},~\pm \sqrt{\frac{2(B\nu - 2A\mu)}{A}},0,0$, that shows it is of type
saddle-saddle-fold.
\item If $B\nu - 2A\mu>0$, and $B\mu - 2C\nu>0$, then $$E_4=(N_1,N_2,M_1,M_2,P_1,P_2) =
(0,0,\sqrt{\frac{B\mu - 2C\nu}{4AC-B^2}},\sqrt{\frac{B\nu - 2A\mu}{4AC-B^2}},0,0)$$ where
the eigenvalues are the roots of quartic polynomial $\alpha x^4 + \beta x^2 +\gamma$, with
\begin{eqnarray*}
\alpha &=& 4AC-B^2>0,\\
\beta &=& 4A(B\mu - 2C\nu)+4C(B\nu - 2A\mu),\\
\gamma &=&  4A(B\mu - 2C\nu)(B\nu - 2A\mu)>0,\\
\end{eqnarray*}
and so if $A>0$ and $C>0$, then $E_4$ is of type hyperbolic-hyperbolic-fold, and if $A<0$
and $C<0$ with $A \,B^{2}+4 A B C +4 A \,C^{2}-2 B^{3}<0$ then $E_4$ is of type
saddle-saddle-fold or center-center-fold, and if $A<0$ and $C<0$ with $A \,
B^{2}+4 A B C +4 A \,C^{2}-2 B^{3}>0$ then $E_4$ is of type hyperbolic-hyperbolic-fold.
\end{itemize}
But if $\nu<0,~\mu>0$, then the equilibria of \eqref{6eqs} are as follows,
\begin{itemize}
\item $E_1=(N_1,N_2,M_1,M_2,P_1,P_2)=(0,0,0,0,0,0)$ with eigenvalues $\pm2\sqrt{\mu}i,~\pm2\sqrt{-\nu},0,0$,
that shows the origin is of type saddle-center-fold.
\item If $C<0$, $E_2=(N_1,N_2,M_1,M_2,P_1,P_2)=(0,0,0,\sqrt{-\frac{\mu}{2C}},0,0)$ then if
$B\mu - 2C\nu>0$, with eigenvalues $\pm \sqrt{2\mu},~\pm \sqrt{-\frac{2(B\mu - 2C\nu)}{C}}i,0,0$,
that shows it is of type saddle-center-fold and if $B\mu - 2C\nu<0$, with eigenvalues
$\pm \sqrt{2\mu},~\pm \sqrt{\frac{2(B\mu - 2C\nu)}{C}},0,0$, that shows it is of type
saddle-saddle-fold.
\item If $A>0$, $E_3=(N_1,N_2,M_1,M_2,P_1,P_2)=(0,0,\sqrt{-\frac{\nu}{2A}},0,0,0)$ then if
$B\nu - 2A\mu>0$, with eigenvalues $\pm \sqrt{-2\nu}i,~\pm \sqrt{\frac{2(B\nu - 2A\mu)}{A}},0,0$,
that shows it is of type saddle-center-fold and if $B\nu - 2A\mu<0$, with eigenvalues
$\pm \sqrt{-2\nu}i,~\pm \sqrt{-\frac{2(B\nu - 2A\mu)}{A}}i,0,0$, that shows it is of
type center-center-fold.
\item If $B\nu - 2A\mu>0$, and $B\mu - 2C\nu>0$, then $$E_4=(N_1,N_2,M_1,M_2,P_1,P_2)=
(0,0,\sqrt{\frac{B\mu - 2C\nu}{4AC-B^2}},\sqrt{\frac{B\nu - 2A\mu}{4AC-B^2}},0,0)$$ where
the eigenvalues are the roots of quartic polynomial $\alpha x^4 + \beta x^2 +\gamma$, with
\begin{eqnarray*}
\alpha &=& 4AC-B^2>0,\\
\beta &=& 4A(B\mu - 2C\nu)+4C(B\nu - 2A\mu),\\
\gamma &=&  4A(B\mu - 2C\nu)(B\nu - 2A\mu)>0,\\
\end{eqnarray*}
and so if $A>0$ and $C>0$, then $E_4$ is of type hyperbolic-hyperbolic-fold,
and if $A<0$ and $C<0$ with $A \,B^{2}+4 A B C +4 A \,C^{2}-2 B^{3}<0$ then $E_4$ is
of type saddle-saddle-fold or center-center-fold, and if $A<0$ and $C<0$ with
$A \,B^{2}+4 A B C +4 A \,C^{2}-2 B^{3}>0$ then $E_4$ is of type hyperbolic-hyperbolic-fold.
\end{itemize}

\subsection{The center-center case for origin}

This case deals with $\nu>0,~\mu>0$, and the equilibria of \eqref{6eqs} are as follows,
\begin{itemize}
\item $E_1=((N_1,N_2,M_1,M_2,P_1,P_2)=(0,0,0,0,0,0)$ with eigenvalues
$\pm2\sqrt{\mu}i,~\pm2\sqrt{\nu}i,0,0$, that shows the origin is of type center-center-fold.
\item If $C<0$, $E_2=(N_1,N_2,M_1,M_2,P_1,P_2)=(0,0,0,\sqrt{-\frac{\mu}{2C}},0,0)$ then if
$B\mu - 2C\nu>0$, with eigenvalues $\pm \sqrt{-2\mu}i,~\pm \sqrt{-\frac{2(B\mu - 2C\nu)}
{C}}i,0,0$, that shows it is of type center-center-fold and if $B\mu - 2C\nu<0$, with
eigenvalues $\pm \sqrt{-2\mu}i,~\pm \sqrt{\frac{2(B\mu - 2C\nu)}{C}},0,0$, that shows it is
of type saddle-center-fold.
\item If $A<0$, $E_3=(N_1,N_2,M_1,M_2,P_1,P_2)=(0,0,\sqrt{-\frac{\nu}{2A}},0,0,0)$
then if $B\nu - 2A\mu>0$, with eigenvalues $\pm \sqrt{2\nu},~\pm \sqrt{-\frac{2(B\nu - 2A\mu)}
{A}}i,0,0$, that shows it is of type saddle-center-fold and if $B\nu - 2A\mu<0$, with
eigenvalues $\pm \sqrt{2\nu},~\pm \sqrt{\frac{2(2\nu - 2A\mu)}{A}},0,0$, that shows it is
of type saddle-saddle-fold.
\item If $B\nu - 2A\mu>0$, and $B\mu - 2C\nu>0$, then $$E_4=(N_1,N_2,M_1,M_2,P_1,P_2)=
(0,0,\sqrt{\frac{B\mu - 2C\nu}{4AC-B^2}},\sqrt{\frac{B\nu - 2A\mu}{4AC-B^2}},0,0)$$ where
the eigenvalues are the roots of quartic polynomial $\alpha x^4 + \beta x^2 +\gamma$, with
\begin{eqnarray*}
\alpha &=& 4AC-B^2>0,\\
\beta &=& 4A(B\mu - 2C\nu)+4C(B\nu - 2A\mu),\\
\gamma &=&  4A(B\mu - 2C\nu)(B\nu - 2A\mu)>0,\\
\end{eqnarray*}
and so if $A>0$ and $C>0$, then $E_4$ is of type hyperbolic-hyperbolic-fold, and if $A<0$
and $C<0$ with $A \,B^{2}+4 A B C +4 A \,C^{2}-2 B^{3}<0$ then $E_4$ is of type
saddle-saddle-fold or center-center-fold, and if $A<0$ and $C<0$ with
$A \,B^{2}+4 A B C +4 A \,C^{2}-2 B^{3}>0$ then $E_4$ is of type hyperbolic-hyperbolic-fold.
\end{itemize}

\subsection{Bifurcation curves}
\noindent

\begin{figure}
\begin{center}
\begin{picture}(210,220)
   \put(0,10){\includegraphics[width=200pt]{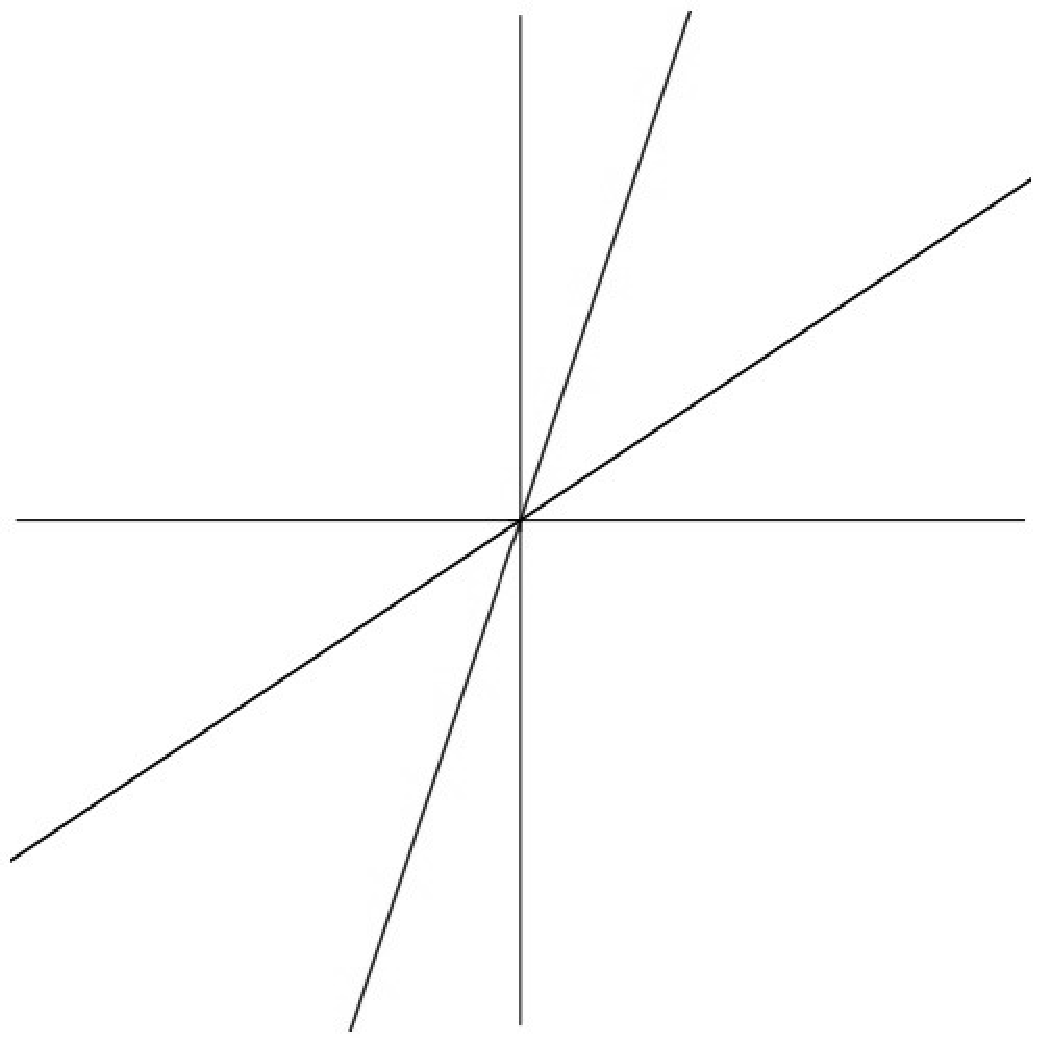}}
   \put(204,105){$\nu$}
   \put(204,170){${\scriptstyle l_1}$}
   \put(135,216){${\scriptstyle l_2}$}
   \put(150,125){${\scriptstyle E_1(cc)}$}
   \put(125,155){${\scriptstyle E_4(hh)}$}
   \put(135,165){${\scriptstyle E_1(cc)}$}
   \put(102,200){${\scriptstyle E_1(cc)}$}
   \put(42,155){${\scriptstyle E_3(cc)}$}
   \put(32,165){${\scriptstyle E_1(sc)}$}
   \put(32,95){${\scriptstyle E_3(cc)}$}
   \put(22,85){${\scriptstyle E_2(sc)}$}
   \put(12,75){${\scriptstyle E_1(ss)}$}
   \put(42,55){${\scriptstyle E_4(hh)}$}
   \put(32,45){${\scriptstyle E_3(cc)}$}
   \put(22,35){${\scriptstyle E_2(sc)}$}
   \put(12,25){${\scriptstyle E_1(ss)}$}
   \put(76,35){${\scriptstyle E_3(cc)}$}
   \put(74,25){${\scriptstyle E_2(cc)}$}
   \put(72,15){${\scriptstyle E_1(ss)}$}
   \put(125,65){${\scriptstyle E_3(sc)}$}
   \put(135,55){${\scriptstyle E_2(cc)}$}
   \put(145,45){${\scriptstyle E_1(sc)}$}
   \put(96,216){$\mu$}
\end{picture}
\end{center}
\caption{
Bifurcations in the case $A>0$, $C>0$, and $B>0$.
\label{BifsACpBp}}
\end{figure}

\begin{figure}
\begin{center}
\begin{picture}(210,220)
   \put(0,10){\includegraphics[width=200pt]{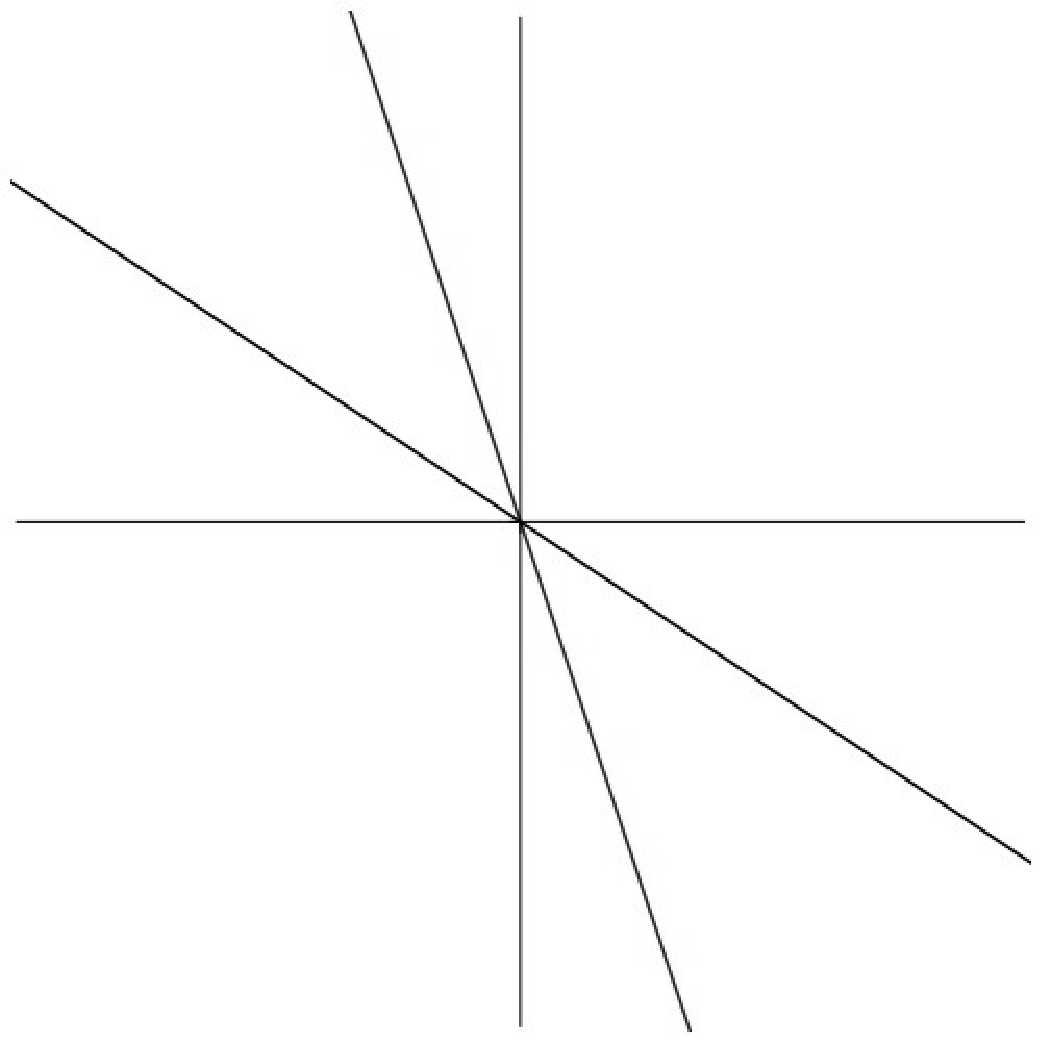}}
   \put(204,105){$\nu$}
   \put(64,214){${\scriptstyle l_1}$}
   \put(0,180){${\scriptstyle l_2}$}
   \put(135,165){${\scriptstyle E_1(cc)}$}
   \put(72,200){${\scriptstyle E_1(sc)}$}
   \put(42,165){${\scriptstyle E_3(sc)}$}
   \put(32,175){${\scriptstyle E_1(sc)}$}
   \put(32,120){${\scriptstyle E_4(hh)}$}
   \put(22,130){${\scriptstyle E_3(sc)}$}
   \put(12,140){${\scriptstyle E_1(sc)}$}
   \put(12,35){${\scriptstyle E_1(ss)}$}
   \put(22,45){${\scriptstyle E_2(sc)}$}
   \put(32,55){${\scriptstyle E_3(sc)}$}
   \put(42,65){${\scriptstyle E_4(hh)}$}
   \put(102,20){${\scriptstyle E_1(sc)}$}
   \put(140,45){${\scriptstyle E_1(sc)}$}
   \put(165,80){${\scriptstyle E_1(sc)}$}
   \put(155,90){${\scriptstyle E_2(sc)}$}
   \put(96,216){$\mu$}
\end{picture}
\end{center}
\caption{
Bifurcations in the case $A>0$, $C>0$, and $B<0$.
\label{BifsACpBn}}
\end{figure}

\begin{figure}
\begin{center}
\begin{picture}(210,220)
   \put(0,10){\includegraphics[width=200pt]{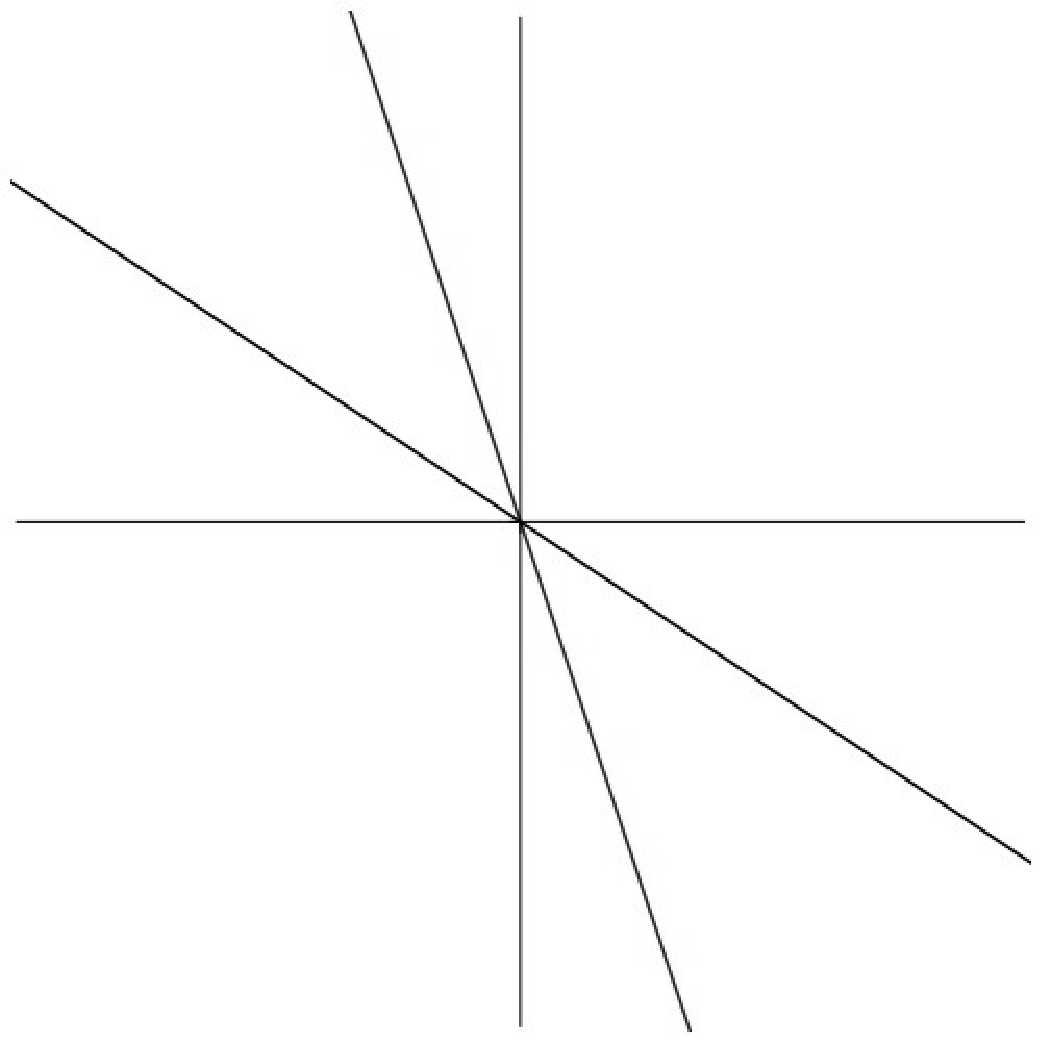}}
   \put(204,105){$\nu$}
   \put(64,214){${\scriptstyle l_1}$}
   \put(0,180){${\scriptstyle l_2}$}
   \put(145,165){${\scriptstyle E_1(ss)}$}
   \put(135,155){${\scriptstyle E_2(cc)}$}
   \put(125,145){${\scriptstyle E_3(sc)}$}
   \put(115,135){${\scriptstyle E_4}$}
   \put(72,205){${\scriptstyle E_1(sc)}$}
   \put(73,197){${\scriptstyle E_2(sc)}$}
   \put(73,190){${\scriptstyle E_4}$}
   \put(52,155){${\scriptstyle E_4}$}
   \put(42,165){${\scriptstyle E_2(sc)}$}
   \put(32,175){${\scriptstyle E_1(sc)}$}
   \put(22,130){${\scriptstyle E_2(ss)}$}
   \put(12,140){${\scriptstyle E_1(sc)}$}
   \put(45,55){${\scriptstyle E_1(ss)}$}
   \put(102,20){${\scriptstyle E_1(sc)}$}
   \put(102,30){${\scriptstyle E_2(cc)}$}
   \put(140,45){${\scriptstyle E_1(sc)}$}
   \put(130,55){${\scriptstyle E_3(sc)}$}
   \put(120,65){${\scriptstyle E_4}$}
   \put(140,45){${\scriptstyle E_1(sc)}$}
   \put(150,95){${\scriptstyle E_4}$}
   \put(160,85){${\scriptstyle E_3(sc)}$}
   \put(170,75){${\scriptstyle E_1(sc)}$}
   \put(96,216){$\mu$}
\end{picture}
\end{center}
\caption{
Bifurcations in the case $A<0$, $C<0$, and $B>0$. If $A \,B^{2}+4 A B C +4 A \,
C^{2}-2 B^{3}<0$ then $E_4$ is of type saddle-saddle or center-center, and if
$A \,B^{2}+4 A B C +4 A \,C^{2}-2 B^{3}>0$ then $E_4$ is of type hyperbolic-hyperbolic-fold.
\label{BifsACnBp}}
\end{figure}

\begin{figure}
\begin{center}
\begin{picture}(210,220)
   \put(0,10){\includegraphics[width=200pt]{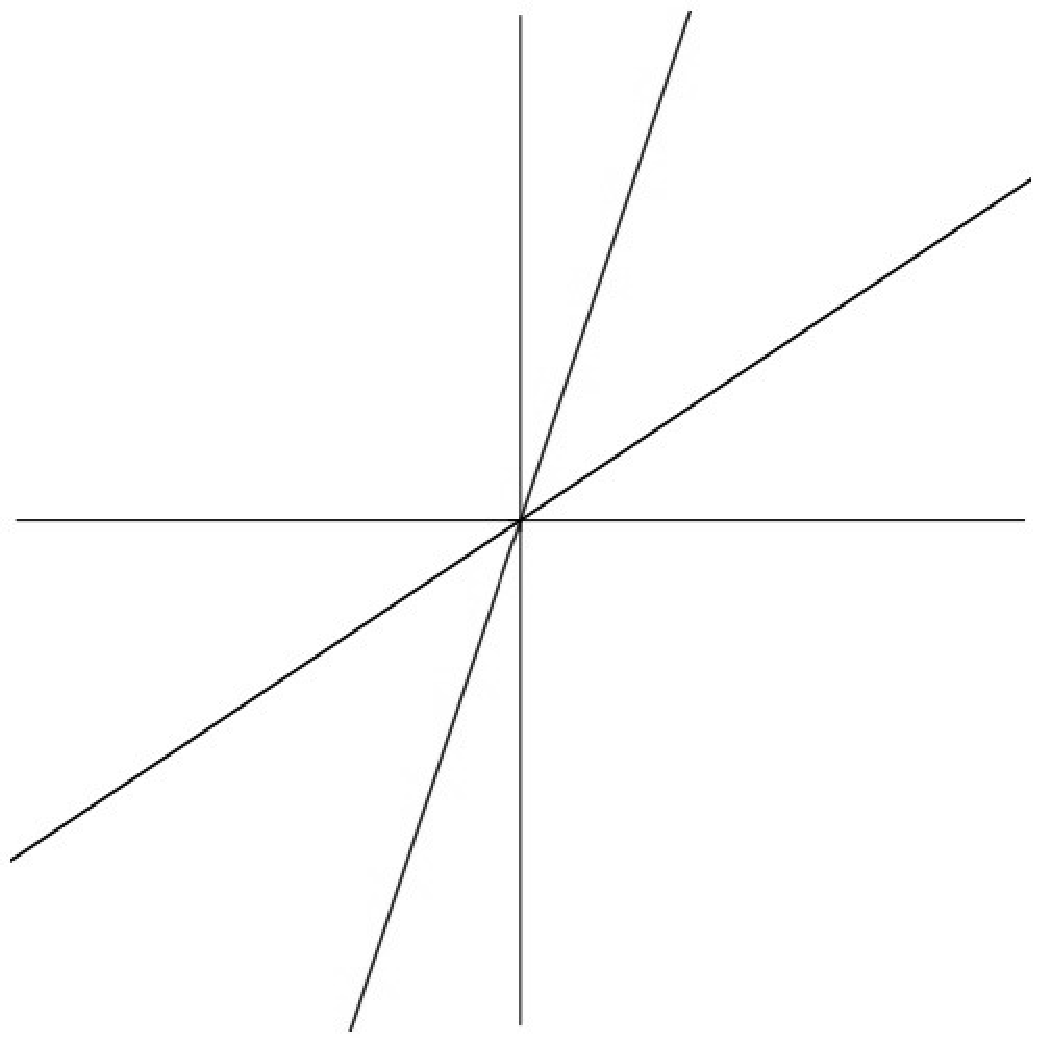}}
   \put(204,105){$\nu$}
   \put(204,170){${\scriptstyle l_1}$}
   \put(135,216){${\scriptstyle l_2}$}
   \put(155,125){${\scriptstyle E_3(ss)}$}
   \put(165,135){${\scriptstyle E_2(cc)}$}
   \put(175,145){${\scriptstyle E_1(cc)}$}
   \put(115,145){${\scriptstyle E_3(sc)}$}
   \put(125,155){${\scriptstyle E_2(cc)}$}
   \put(135,165){${\scriptstyle E_1(cc)}$}
   \put(102,203){${\scriptstyle E_1(cc)}$}
   \put(101,197){${\scriptstyle E_2(sc)}$}
   \put(100,191){${\scriptstyle E_3(sc)}$}
   \put(99,185){${\scriptstyle E_4}$}
   \put(42,155){${\scriptstyle E_2(sc)}$}
   \put(32,165){${\scriptstyle E_1(sc)}$}
   \put(22,85){${\scriptstyle E_1(ss)}$}
   \put(52,65){${\scriptstyle E_4}$}
   \put(42,55){${\scriptstyle E_1(ss)}$}
   \put(72,15){${\scriptstyle E_1(ss)}$}
   \put(125,65){${\scriptstyle E_3(ss)}$}
   \put(135,55){${\scriptstyle E_1(sc)}$}
   \put(96,216){$\mu$}
\end{picture}
\end{center}
\caption{
Bifurcations in the case $A<0$, $C<0$, and $B<0$. If $A \,B^{2}+4 A B C +4 A \,C^{2}-2 B^{3}<0$
then $E_4$ is of type saddle-saddle or center-center, and if $A \,B^{2}+4 A B C +4 A \,
C^{2}-2 B^{3}>0$ then $E_4$ is of type hyperbolic-hyperbolic-fold.
\label{BifsACnBn}}
\end{figure}

Since we have the integral
$$I= N_1 + N_2 +\nu M_1 + \mu M_2 + A M_1^2+ B M_1M_2+ C M_2^2,$$
the dynamics can be reduced on an invariant $4$-torus. By this invariant manifold,
we can consider the reduced system with non-zero eigenvalues. Since $A$ and $C$ have
the same signs, for the bifurcations, we can consider the following cases,
\begin{enumerate}
\item[case (1)] $A>0$, $C>0$, and $B>0$,
\item[case (2)] $A>0$, $C>0$, and $B<0$,
\item[case (3)] $A<0$, $C<0$, and $B>0$,
\item[case (4)] $A<0$, $C<0$, and $B<0$.
\end{enumerate}
Consider the following lines,
$$l_1=\{(\nu,\mu)~:~B\nu - 2A\mu=0\},$$
$$l_2=\{(\nu,\mu)~:~B\mu - 2C\nu=0\}\enspace.$$
For each case, we have the bifurcation diagram on $(\nu,\mu)$-plane for the reduced system
in Figs. \ref{BifsACpBp}, \ref{BifsACpBn}, \ref{BifsACnBp}, and \ref{BifsACnBn}, respectively.
Let we mean $ss$, $sc$, $cc$ and $hh$, saddle-saddle, saddle-center, center-center and
hyperbolic-hyperbolic respectively. In the figures $E_j(\ast\dag)$ show the equilibrium $E_j$
of type $\ast\dag$ for $j=1,2,3,4$. We notice that in all cases when $\nu<0,~\mu<0$, the origin
is of type saddle-saddle (i.e. $ss$).

\noindent

\noindent

\noindent

\newpage

\section{Nonintegrability of the truncated system as $\zeta_1=\zeta_2=0$}

In this section we study the Hamiltonian system with the truncated Hamiltonian of
the fourth degree
$$\begin{array}{l}
H= \omega_1(p_2q_1 - p_1q_2)+\frac12(q_1^2 + q_2^2) +
\omega_2(p_4q_3 - p_3q_4)+ \frac12(q_3^2 + q_4^2) +
\frac{\eps_1}{2}(p_1^2 + p_2^2) + \frac{\eps_2}{2}(p_3^2 + p_4^2)+\\
\frac 14[A(p_1^2 + p_2^2)^2 + 2B(p_1^2 + p_2^2)(p_3^2 + p_4^2)+
C(p_3^2 + p_4^2)^2].
\end{array}
$$
As we discussed above, for this Hamiltonian system functions $S_1 = p_2q_1 - p_1q_2$ and
$S_2 = p_4q_3 - p_3q_4$ are integrals with their related Hamiltonian systems
defining the actions of the group $S^1$. So they generate the action of 2-torus in the
phase space $\R^8$. The reduction w.r.t. this action leads to a 2-DOF Hamiltonian system
\cite{Smale,MW}. The first question here: is this system integrable, what is its structure?
Notice that this system with four degrees of freedom has three almost everywhere independent
involutive integrals: $H, S_1, S_2$. For its integrability in the Liouville sense one needs
to have one more integral being in involution with these three ones. But we shall
show that generically no more integrals exist for the system.

The reduction of the system is performed most easily in new symplectic coordinates,
similar to \cite{sok}
\begin{align}
p_1 &= r_1\cos\varphi_1,&q_1 &= Q_1\cos\varphi_1 + \frac{k_1}{r_1}\sin\varphi_1,\\
p_2 &= r_1\sin\varphi_1,&q_2 &= Q_1\sin\varphi_1 - \frac{k_1}{r_1}\cos\varphi_1,\\
p_3 &= r_2\cos\varphi_2,&q_3 &= Q_2\cos\varphi_2 + \frac{k_2}{r_2}\sin\varphi_2,\\
p_4 &= r_2\cos\varphi_2,&q_4 &= Q_2\sin\varphi_2 - \frac{k_2}{r_2}\cos\varphi_2,
\end{align}
in which the symplectic 2-form looks as follows
$$
dq_1\wedge dp_1 + dq_2\wedge dp_2 + dq_3\wedge dp_3 + dq_4\wedge dp_4 =
dQ_1\wedge dr_1 + dQ_2\wedge dr_2 + d\varphi_1\wedge dk_1 + d\varphi_2\wedge dk_2.
$$
These variables are expressed via invariants $(M_i,N_i,P_i,S_i)$. Indeed,
we have
\begin{equation}\label{inv}
\begin{array}{l}
k_1=S_1,\;k_2=S_2,\;r_1Q_1 = P_1,\;r_2Q_2 = P_2,\;r_1^2 = 2M_1,\;r_2^2 = 2M_2,\\
Q_1^2+k_1^2/r_1^2 = 2N_1,\;Q_2^2+k_2^2/r_2^2 = 2N_2,\;
\end{array}
\end{equation}
that implies the syzygies $4M_iN_i-P_i^2 = S_i^2,$ $i=1,2.$
\begin{remark}
Simply, we can notice that the parameters $k_1$ and $k_2$ are same as $\zeta_1$ and
$\zeta_2$, respectively, in the previous section.
\end{remark}

Despite the functions defining the change of variables have singularities at points of
two planes $p_1=p_2=0$ and $p_3=p_4=0$, these variables are very convenient. The Hamiltonian
$H$ in these variables becomes
\begin{equation}\label{sk}
H = \omega_1k_1 +\omega_2k_2 + \frac12(Q_1^2 + k_1^2/r_1^2) + \frac12(Q_2^2 + k_2^2/r_2^2) +
\frac{\eps_1}{2}r_1^2 + \frac{\eps_2}{2}r_2^2 + \frac 14[Ar_1^4 + 2Br_1^2r_2^2 +
Cr_2^4].
\end{equation}

\begin{remark}
Simply, we can notice that the parameters $\eps_1$ and $\eps_2$ are related to $\nu$ and
$\mu$ respectively in the previous section.
\end{remark}

So, $k_1, k_2$ become new parameters (they are integrals), $\varphi_1, \varphi_2$ are cyclic
variables, whence the system casts as the 2-DOF system in variables $(Q_1, Q_2, r_1, r_2)$ which
depends on four parameters $k_1,k_2, \eps_1, \eps_2$. One needs to study this system and
understand how its structure varies with parameters. We consider below the values of
$r_1,r_2$ taking any sign, not only positive, ignoring so far its geometric
sense.

The first task here is to investigate the regular case $k_1 = k_2 =0,$ that corresponds to
zeroth level of integrals $S_1, S_2$ and to the reduced system on the 6-cone. Namely, in each
4-subspace of the phase space $\R^8$ the related quadratic integral defines a 3-dimensional
cone being topologically a cone over 2-torus. The 6-dimensional set obtained
in $\R^8$ is invariant w.r.t. the action of the group $\R^2$ generated by commuting Hamiltonian
vector fields $X_{S_1}, X_{S_2}$. The orbits of this $\R^2$-action are mostly two-dimensional
tori, but they degenerate on some subsets into the closed orbits. Reduction of the 6-dimensional
set w.r.t. these two integrals leads to the 4-dimensional symplectic manifolds with the reduced
flow \cite{Smale,MW}. The most simply this scheme works for the case $k_1 = k_2 =0.$ Then the
reduced space in these coordinates is merely $\R^4.$

The Hamiltonian of the system at $k_1 = k_2 =0$ is a quartic polynomial
\begin{equation}\label{zero}
H_0 = \frac12(Q_1^2 + Q_2^2) + \frac{\eps_1}{2}r_1^2 + \frac{\eps_2}{2}r_2^2 +
\frac 14[Ar_1^4 + 2Br_1^2r_2^2 + Cr_2^4]
\end{equation}
and the problem of its integrability was partially investigated in several papers
\cite{MS,LS,Yo,Ziglin}. There it was discovered that for some isolated points in the parameter
space integrability indeed occurs and additional integrals were found for these cases.
Nonetheless, in general these systems are not integrable, as we shall see.
\begin{example}
For the case $\eps_1=\eps_2 = \eps,$ $A=C, B=3A$ both integrals are as
follows \cite{LS}
$$
H = \frac12(Q_1^2+Q_2^2) + \frac{\eps}{2}(r_1^2+r_2^2) +
\frac{A}{4}(r_1^4+6r_1^2r_2^2+r_2^4),\;
K = Q_1Q_2 + \eps r_1r_2 + Ar_1r_2(r_1^2+r_2^2).
$$
For the case $A = B = C$, i.e. $H = \frac12(Q_1^2+Q_2^2) + \frac{\eps_1}{2}r_1^2 +
\frac{\eps_1}{2}r_2^2 + \frac{A}{4}(r_1^2+r_2^2)^2$, one has the Garnier system \cite{Gar}
(see, for instance \cite{Per}), the second integral is \cite{LS}
$$
K = \frac{\eps_2-\eps_1}{2}Q_1^2 + \eps_1\frac{\eps_2-\eps_1}{2}r_1^2 +
\frac{A}{4}[(\eps_2-\eps_1)r_1^2(r_1^2+r_2^2) + (r_1Q_2-r_2Q_1)^2]
$$
There are two more integrable cases with fifth and eighth order
polynomials as additional integrals, see \cite{LS,LS1}.
\end{example}

Hamiltonian system (\ref{zero}) possesses discrete symmetries. It is reversible
w.r.t. two involutions
$$\begin{array}{l}
L_1: (Q_1,Q_2,r_1,r_2)\to (-Q_1,Q_2,r_1,r_2),\\
L_2: (Q_1,Q_2,r_1,r_2)\to (Q_1,-Q_2,r_1,r_2).
\end{array}
$$
Besides, it is symmetric w.r.t. two involutions: $R_1: (Q_1,Q_2,r_1,r_2)\to
(-Q_1,Q_2,-r_1,r_2)$ and $R_2: (Q_1,Q_2,r_1,r_2)\to (Q_1,-Q_2,r_1,-r_2).$
The system has two invariant 2-planes $Q_2=r_2 = 0$ and $Q_1=r_1 =
0,$ on which the restriction of the system gives the phase portrait of the
Duffing equation \cite{Morozov}. In fact, these systems depending on parameters $(\eps_1,\eps_2)$ were
investigated in the previous section. This study suggests that homoclinic orbits of
the equilibrium $O$ at the origin can exist for the corresponding signs of the ratios
$\eps_1/A$ and $\eps_2/C.$

The Hamiltonian system with the Hamiltonian $H_0$ takes the form of two coupled anharmonic
oscillators
\begin{align}\label{sys_0}
\dot r_1 &= Q_1,&\dot Q_1 &= -\eps_1r_1 - Ar_1^3 - Br_1r_2^2,\nonumber\\
\dot r_2 &= Q_2,&\dot Q_2 &= -\eps_2r_2 - Br_1^2r_2 - Cr_2^3.
\end{align}
The system satisfies all assumptions of the paper \cite{BLT}, if both $\eps_1,\eps_2$
are negative and $A, C$ are positive, so all its conclusions hold.
As particular structures of this system, two homoclinic orbits forming figure
``eight'' belong to each invariant 2-planes. They are given by formulas
\begin{equation}\label{hom}
r_1 = \sqrt{-\frac{2\eps_1}{A}}\frac{1}{\cosh(\sqrt{-\eps_1}t)},\;
r_2 = \sqrt{-\frac{2\eps_2}{C}}\frac{1}{\cosh(\sqrt{-\eps_2}t)}.
\end{equation}
We shall use this later.

Let us start with the investigation of the type of the equilibrium at the
origin, how it changes when parameters $\eps_1, \eps_2$ vary near zero
values. The quadratic part of the Hamiltonian $H_0$ is the sum of two
uncoupled partial Hamiltonians $H_0^{(2)} = \frac 12 (Q_1^2
+\eps_1r_1^2)+$ $\frac 12 (Q_2^2 +\eps_2r_2^2)$. So, we get the direct
product of two saddles on each 2-plane $(Q_1,r_1),$ $(Q_2,r_2),$ if both $\eps_1,
\eps_2$ are negative, of a saddle and a center, if the product
$\eps_1\eps_2$ is negative, and the direct product of centers on each 2-plane
$(Q_1,r_1),$ $(Q_2,r_2),$ if both $\eps_1,\eps_2$ are positive. In $\R^4$ we get in
the linear approximation respectively a saddle equilibrium, a saddle-center and
an elliptic point.

First we assume the quadratic form $A\xi_1^2 + 2B\xi_1\xi_2 + C\xi_2^2$ be positive or
negative definite, i.e. $AC-B^2>0$. Then the levels of the Hamiltonian $H=c$ for $c>0$ large
enough (or for negative $c$ large enough in modulus) are compact submanifolds diffeomorphic
to the 3-sphere $S^3.$ Thus, in this case all orbits of the system (\ref{sys_0}) in these
levels are continued onto the whole $\R$ and their limit sets are compact. This
assumption, in particular, implies that constants $A,C$ has the same sign, if
$B\ne 0$: $AC>B^2 \ge 0.$ We always assume $B\ne 0$, otherwise the
subsystems are decoupled.

\subsection{The saddle-center case}

\begin{figure}[h]
	\centering
	\includegraphics[width=0.6\linewidth]{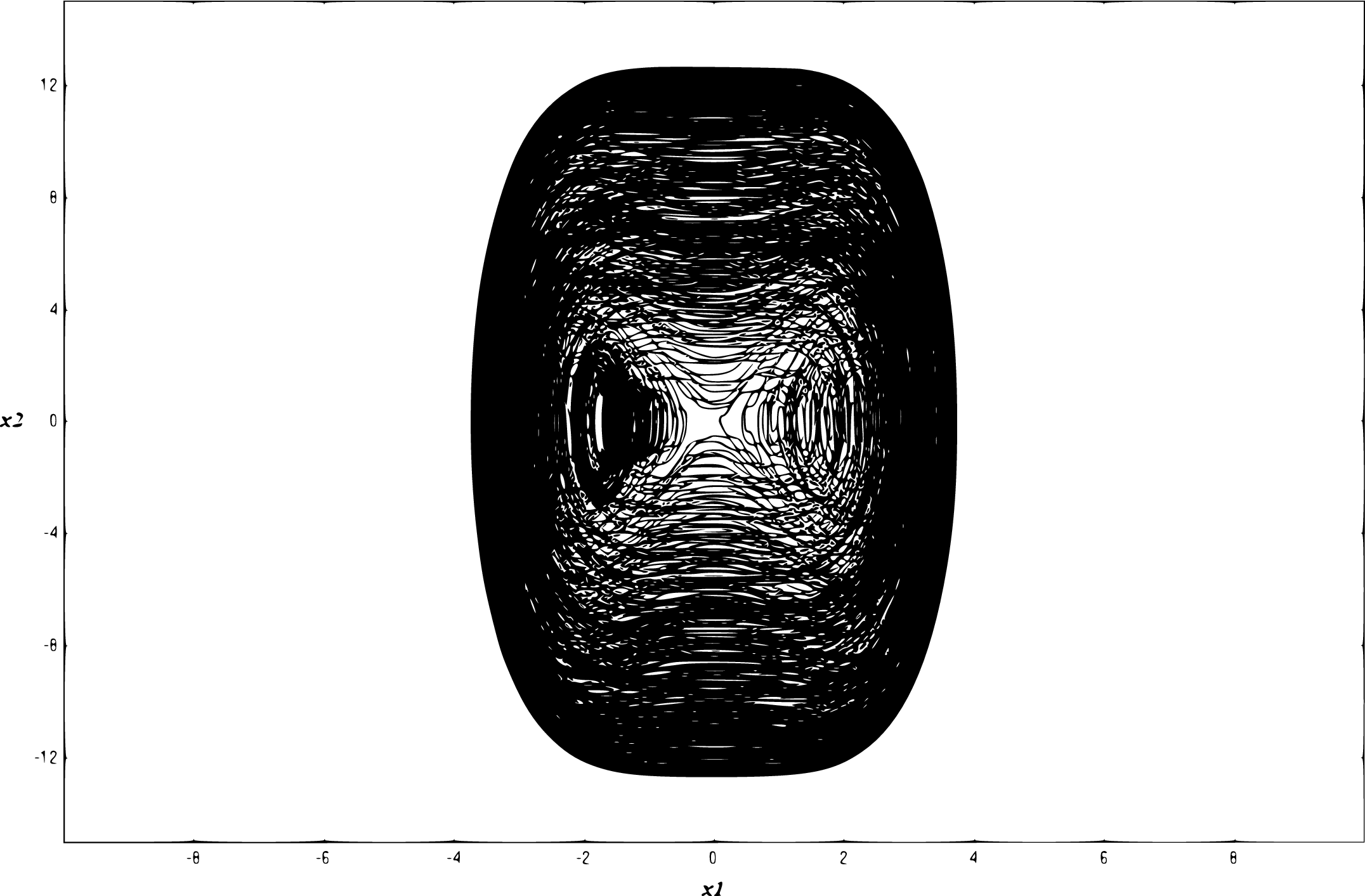}
	\caption{Projection on the plane $r_1,r_2.$}
\label{fig:7}
\end{figure}

For the case of a saddle-center at the origin the system (\ref{sys_0}) possesses
a ``figure eight'' on one invariant plane and a center at the origin on another invariant
plane (in this case, if the ratio $-\eps_1/A$ or $-\eps_2/C$, respectively, is
positive, the related invariant 2-plane have also two heteroclinic orbits joining
two saddles). The orbit behavior here can be very complicated that is seen on the plot
(see Fig.~\ref{fig:7}) obtained by the numerical integration of the system (\ref{zero})
at the parameters
$\eps_1=-0.3,\eps_2=0.3, A=2.B=-1, C=3.$

For this case we are able to prove rigorously that the system generically is not integrable
and therefore the complicated behavior is a genuine feature of this system. To this end, we
recall results from \cite{Le,Grotta}. The orbit structure of a Hamiltonian system near
a homoclinic orbit of a saddle-center was studied first in \cite{Le} (without any assumption
that the homoclinic orbit lies on an invariant submanifold), where it was proved
that if some genericity condition of the linearized system at the homoclinic solution holds,
then any Lyapunov's small periodic orbit $\gamma_c$ on the local center manifold possesses four
transverse homoclinic orbits and therefore the system is not integrable (due to results by
Cushman \cite{Cu} the system cannot have an analytic integral, see also other methods to get
this result \cite{Ziglin,Morales-Ruiz}). The important consequence of this
result is the fact that non-integrability can be detected by investigating
the structure of the linearized system on the homoclinic orbit only.

For a Hamiltonian system with the homoclinic orbit on the invariant plane, like (\ref{sys_0}),
the genericity condition can be verified directly using the linearized system on
the homoclinic solution. This is very important that non-integrability of a nonlinear system
can be guaranteed by the structure related linearized system on the homoclinic solution. More
results in this approach can be found in papers \cite{Ziglin,Morales-Ruiz,MR_book,Yagasaki}.
In the system (\ref{sys_0}) the linearized system splits due to presence of the invariant
plane and the corresponding normal linear subsystem is defined. The scattering map for
this subsystem answers on the question of integrability \cite{Grotta}. In more details,
consider an analytic Hamiltonian system with two degrees of freedom
\begin{equation}\label{GR}
\dot x = H_y,\; \dot y = -H_x,\; x=(x_1,x_2),\; y=(y_1,y_2),
\end{equation}
with a classical analytic Hamiltonian of the form $H = (y_1^2 + y_2^2)/2 + V(x_1,x_2)$,
where the potential $V$ satisfies the conditions:
\begin{enumerate}
\item $V = \frac12[-\nu^2 x_2^2 + \omega^2 x_2^2]+ O((x_1^2 +
x_2^2)^{3/2})$ with real nonzero $\nu,\omega$;
\item $\partial_{x_2}V(x_1,0)\equiv 0$ for any $x_1\in \R$;
\item the equation $V(x_1,0)=0$ has a nontrivial non-degenerate solution
$x_1(c)$ and no solutions on the interval $(0,x_1(c)).$
\end{enumerate}
Such Hamiltonian system has the following properties:

- the origin is an equilibrium point of saddle-center type, namely it is associated
with a pair of real $\pm\nu$ and a pair of imaginary eigenvalues $\pm i\omega$;

- there is a homoclinic orbit to the origin (it belongs to the invariant plane
$x_2 = y_2 = 0$).

For this system the following theorem holds
\begin{theorem}[Grotta Ragazzo]
Let $\Gamma:\R \to \R$ be the $x_1$-component of the solution to the system (\ref{GR})
being a homoclinic orbit of the equilibrium at the origin. If the ``reflection coefficient''
$R$ associated to the scattering problem of the linear differential system:
$$\ddot z = -[\partial_{x_2x_2}V(\Gamma,0)]z
$$
is different from zero, then the system is non-integrable.
\end{theorem}
Non-integrability here means the same as in \cite{Le,Grotta,MHR}, i.e. the existence of four
transversal homoclinic orbits to each Lyapunov' periodic orbit.

The calculation of the scattering coefficient $R$ is a rather complicated
problem, generally speaking, but it can be resolved for the particular type of the potential
$V$. Namely, suppose $V$ satisfies:
$$
\begin{array}{l}
V(x_1,0) = -\frac12\nu^2 x_1^2 + \frac{\alpha}{4}x_1^4,\;\alpha > 0,\\
\partial_{x_2x_2}V(x_1,0) = \omega^2 + \beta x_1^2.
\end{array}
$$
If
$$
\frac{\beta}{\alpha}\ne \frac{l(l-1)}{2}, \;l\in \mathbb N,
$$
then the system (\ref{GR}) with this $V$ is non-integrable.

For our case the needed coefficients are $\alpha = A,$ $\beta = B$. So, we have
\begin{cor}
For all ratios $B/A$, except for integers $\frac{l(l-1)}{2}, \;l\in \mathbb N$,
the system (\ref{GR}) is non-integrable.
\end{cor}

An invariant formulation of the genericity condition for a Hamiltonian system with a
homoclinic orbit $\Gamma$ of a saddle-center $p$, independently on whether the homoclinic orbit
of the system lies on an invariant plane or not, was found in \cite{KL} and
its generalization to the case of $n$ degrees of freedom Hamiltonian system with a homoclinic
orbit to a 1-elliptic equilibrium (having a pair of simple imaginary eigenvalues and others
out of the imaginary axis) was found in \cite{KL1}.

\begin{figure}
\includegraphics[scale=0.6]{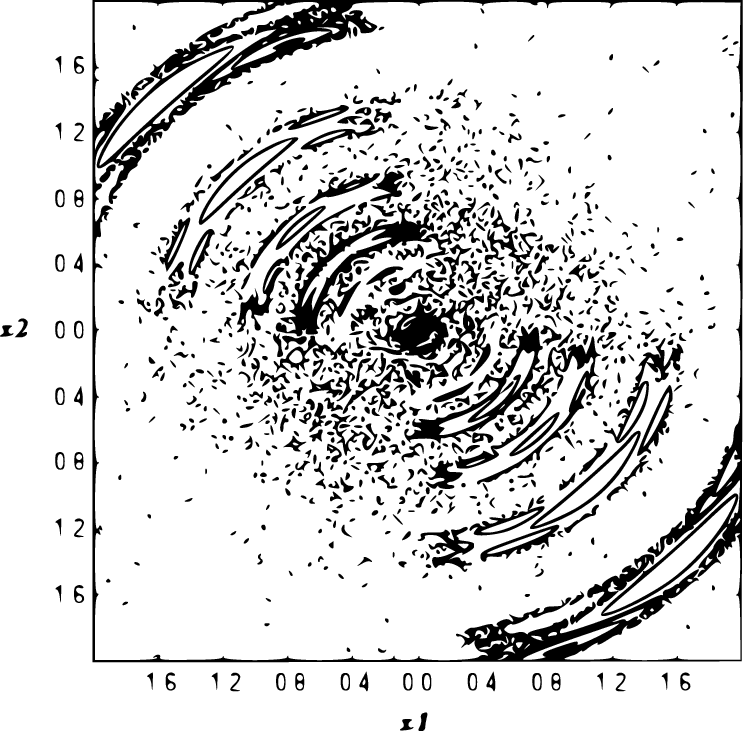} \includegraphics[scale=0.6]{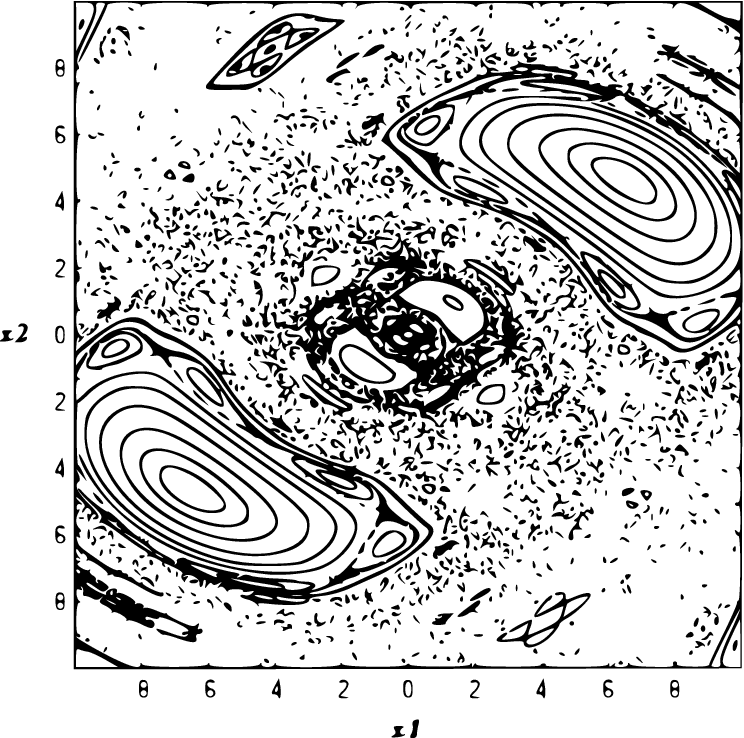}
\caption{Orbit behavior of the Poincar\'e map on a cross-section to a homoclinic
orbit of a saddle-center}
\label{fig:8}
\end{figure}

Except for the existence in a neighborhood of $\Gamma$ of four transverse homoclinic orbits
to each small Lyapunov's periodic orbit on the center manifold of the
saddle-center, there is a rather standard orbit structure on the singular
level $H=H(p)$, namely \cite{KL1,GR2}
\begin{theorem}
The singular level of the Hamiltonian contains four countable families of
periodic orbits accumulating to the homoclinic orbit $\Gamma$. Two of
these families consists of saddle periodic orbits, but the types of two
other families depends on the value of some constant composed of
eigenvalues and coefficients of the linearized global map: they can consist
of elliptic periodic orbits but can of hyperbolic ones.
\end{theorem}
The illustration of such behaviors is presented on the Fig.~\ref{fig:8}
obtained by the simulation of the related Poincar\'e map on the cross-section to $\Gamma$.

\subsection{The saddle case}\label{6.2}

In the case of a saddle at the origin for the system (\ref{zero}) (i.e. when both
$\eps_1, \eps_2$ are negative) the initial system has at the origin the equilibrium of
the double saddle-focus type (see, Sec. \ref{origin}). If $A,C>0$, the system (\ref{zero})
has four homoclinic orbits to the saddle (two ``figures eight'' on each invariant 2-plane)
in the level $H_0 = 0$. Such homoclinic orbit corresponds to the 3-dimensional invariant set
(a ``homoclinic skirt'') of the truncated system in $\R^8$ with the quartic Hamiltonian
(\ref{truncated-normal-form}) filled with homoclinic orbits of the double saddle-focus which
are generated by the action of the torus $\mathbb T^2$ on the points of the homoclinic orbit
(\ref{hom}) in the related plane. For the case of the double saddle-focus the homoclinic skirt
belongs to the intersection of four-dimensional stable $W^s(O)$ and instable $W^u(O)$ manifolds
of $O$, both of them belong to the five-dimensional joint level $H =  S_1 = S_2 = 0.$ Hence,
this intersection can be transverse along this 3-dimensional skirt. This reminds the situation
in integrable Hamiltonian systems with two degrees of freedom which have a saddle equilibrium
(with simple real eigenvalues $\pm\lambda_1, \pm\lambda_2$). There the saddle can have four
transverse homoclinic orbits in the singular level of a Hamiltonian \cite{Deva,LU}. What is
interesting and it was mentioned first in \cite{Deva} that an integrable Hamiltonian system
can have such transverse homoclinic orbits (the known Neumann's integrable Hamiltonian system
was presented in \cite{Deva}). Later this was shown to be a general phenomenon for an integrable
Hamiltonian system with a saddle-saddle singular point being simple in the Vey-Eliasson sense
\cite{LU}. Of course, such situation can be met and in non-integrable systems with the saddle
equilibrium, so such situation cannot be a criterion for non-integrability
of a Hamiltonian system, in contrast with the case of a transverse homoclinic orbit of
a saddle-focus equilibrium \cite{Dev,Lesf,Ler1}.

For the system (\ref{sys_0}) we also expect existing super-homoclinic orbits in
$\R^4$,lying in the critical level of the Hamiltonian $H_0$, and infinitely many multi-round
homoclinic orbits to $O$ for systems in some regions of the parameter plane $(\eps_1,\eps_2)$
\cite{BLT}. Recall \cite{EKTS,TSh} when a symmetric Hamiltonian system has a saddle equilibrium
with the leading eigenvalue being real simple and has several transverse homoclinic
orbits, then as was shown in \cite{TSh} (based on an earlier work \cite{EKTS}) that
symmetries can lead to the emergence of the so-called super-homoclinic orbits which,
in turn, serve as limits of infinite series of multi-pulse homoclinic loops. Namely, as shown
in \cite{TSh}, in certain situations, the homoclinic loops or bunches of homoclinic loops
can have stable or unstable invariant manifolds; the super-homoclinics correspond
to orbits lying in the intersection of these manifolds. To explain this,
notice that two homoclinic loops forming the ``figure eight'' belong to
the invariant plane and this plane in hyperbolic in the transverse
direction due to the saddle equilibrium at the origin in the transverse invariant
plane. Consider 2-dimensional unstable manifold $W^u(O)$ of the saddle in
$\R^4$. It is continued by the flow of (\ref{sys_0}) and because of
existence of the homoclinic loop it returns to a neighborhood of $O$ and
intersects 2-dimensional stable manifold $W^s(O)$. The existence of
invariant planes implies that the homoclinic orbit leaves neighborhood of
$O$ along leading direction and returns along non-leading direction or
vice versa. Here we assume that $|\eps_1|\ne |\eps_2|$, then eigenvalues of
$O$ are different.

In the levels other than zeroth, we have several saddle periodic orbits which belong to
one-parameter families of periodic orbits accumulating each to its own homoclinic orbits
and to the ``figure eight'' in the related invariant 2-plane. Definitely, on the plane with
the figure eight we have for each nonzero level of the Hamiltonian a pair of periodic orbits
lying each inside of the figure eight eyes. These periodic orbits are saddle ones on the
related level of the Hamiltonian $H_0$. We expect these two saddle periodic orbits are
connected by four heteroclinic orbits forming a heteroclinic connection with complicated
dynamics.

In order to show non-integrability of the system for the saddle case, we
follow to the method due to Ziglin \cite{Ziglin} and elaborated further by
Simo, Morales-Ruiz, Ramis and others \cite{Morales-Ruiz}. This is based on
the investigation of the linearized equation at some particular solution of
the system. Here as such solution we choose a periodic solution inside of
a homoclinic loop on the invariant plane $Q_1=r_1=0.$ As is known, these
solutions are expressed via Jacobi elliptic functions $r_2 = g(t,h) = m(h)dn(s(h)t,k)$
\cite{JEL} with positive smooth functions $m(h), s(h)$ depending on the value of energy
level $h$ where the periodic orbit belongs to. Because the plane is invariant, the linearized
system splits and normal variational equation (NVE) is as follows
$$
\ddot r_1 + [\eps_1 + Bm^2(h) - Bm^2(t,h)k^2(h)sn^2(s(h)t,k)]r_1 = 0,
$$
where $sn(\tau,k)$ denote the Jacobi elliptic sinus, $k$ is its modulus. Thus we have
the Lam\'e type normal variational differential equation \cite{BE,MS}. For
this case the following theorem is valid \cite{MS}.
\begin{theorem}
If the NVE along the solution $\Gamma$ is of the Lam\'e type and falls
outside of cases (i), (ii), (iii) , then the Hamiltonian system has not
meromorphic first integrals independent of the Hamiltonian in a
neighborhood of $\Gamma$.
\end{theorem}
The cases (i), (ii), (iii) enumerate several discrete cases when the
Lam\'e equation is solvable in the differential Galois theory sense
(see, for instance, \cite{MS}).

As a signature of the orbit structure that allows to say about non-integrability
of the system a presence of super-homoclinic orbit near a pair of homoclinic orbits
forming a figure eight can be taken. It is well known that the orbit behavior near
an individual homoclinic orbit $\Gamma$ of a saddle is very simple \cite{VF}: there
is only a one-parametric family of saddle periodic orbits accumulating to $\Gamma$,
no other orbits lying entirely in a neighborhood of $\Gamma$ exist.
Indeed, the integrability of a Hamiltonian system with two degrees of freedom (i.e.
an existence of a smooth additional integral independent almost everywhere with the
Hamiltonian) having a saddle equilibrium with transverse homoclinic orbits implies
that all other orbits on the unstable manifold of $\Gamma$ are homoclinic ones \cite{LU}
(i.e. stable and unstable manifolds of the saddle coalesce), if the common connected level
of two integrals, Hamiltonian and an additional one, form a compact set with the only
equilibrium.

The Hamiltonian system (\ref{sys_0}) falls into those studied in
\cite{BLT} where the super-homoclinic orbits where proved to exist. So,
this gives us the supplementary evidence of non-integrability of the
system.

\subsection{The elliptic case}

When the origin is an elliptic equilibrium, i.e. $\eps_1 > 0, \eps_2 > 0,$ then
this equilibrium $O$ is Lyapunov stable in a neighborhood of the order
$\sqrt{\eps_1^2+\eps_2^2}$, since the Hamiltonian $H_0$ (\ref{zero}) is positive
definite near $O$ and levels $H_0 = c > 0$ for small enough $c$ are nested 3-dimensional
spheres. Both invariant 2-planes in case if $A>0, C>0$ are filled with families of periodic
orbits enclosing the center equilibrium -- Lyapunov's families. Locally near the equilibrium
these periodic orbits are elliptic ones, i.e. their essential multipliers
(different from the double unit) consist of two complex conjugate numbers $\exp[\pm i\alpha]$
on the unit circle.. So we expect the KAM theorem \cite{Arnold} works
here and a positive measure sets of two-dimensional tori exist nearby. But this depends also
on the existence of resonances between frequencies $\sqrt{\eps_1}, \sqrt{\eps_2}.$
If some of constants $A,C$ is negative, then the related invariant plane
contains two symmetric saddle equilibria connected by two heteroclinic
orbits forming a heteroclinic connection inside of which all orbits are
periodic but orbits out of this connection leave to infinity in the
finite time.

The local orbit structure described can exist for both integrable and
nonintegrable systems, the difference between them consists in the more
subtle details like a behavior near resonant torus with rational rotation
numbers, i.e. in the so-called resonant zones \cite{Morozov}. There in a nonintegrable
system usually exist saddle periodic orbits whose stable and unstable
two-dimensional manifolds intersect each other forming homoclinic tangles
with complicated orbit behavior. For the case we study the
nonintegrability can be also detected by means of tools of differential
Galois theory \cite{Morales-Ruiz} based on periodic orbits of two
Lyapunov's families like in the Sec.\ref{6.2}.

Each Lyapunov's periodic orbit, due to the action of the torus $T^2$, gives rise to a family
of three-dimensional invariant tori in $\R^8$ for the truncated Hamiltonian of the fourth order.

\section{Supplement}

In this Section we discuss the 2-mode system derived in \cite{KuLe} and present more detailed
study of bifurcations related with the birth of small homoclinic orbits of the equilibrium at
the origin.  The system is the following
\begin{equation}\label{8}
\begin{array}{l}
q'_1 = q_2,\\
q'_2 = p_2 - q_1,\\
q'_3 = q_4,\\
q'_4 = p_4- (1-\omega^2)q_3,\\
p'_1 = p_2-\alpha q_1 - \frac{\beta}{\sqrt{2}}(q_1^2+q_3^2) +
\frac12q_1^3 + \frac32 q_1q_3^2,\\
p'_2 = - p_1,\\
p'_3 = (1-\omega^2)p_4 - \alpha q_3 - \beta\sqrt{2}q_1q_3 + \frac32 q_3(q_1^2 +
\frac12 q_3^2),\\
p'_4 = -p_3 .
\end{array}
\end{equation}
It is a Hamiltonian system with the Hamiltonian
$$
\begin{array}{l}
H = p_1q_2-p_2q_1 + [p_3q_4-(1-\omega^2)p_4q_3]+
\frac12(p_2^2+p_4^2)+\frac{\alpha}{2}(q_1^2+q_3^2)+\\
\hspace{1.cm}\frac{\beta}{\sqrt{2}}[\frac13 q_1^3 + q_1q_3^2] -
\dst{\frac{1}{16}[2q_1^4 + 3q_3^4 + 12q_1^2q_3^2]}.
\end{array}
$$
For the equilibrium $O$ at the origin the matrix of its linearized system as $\alpha = 0$ has
two double pure imaginary eigenvalues $\pm i,$ $\pm i \sqrt{1-\omega^2}$ when $|\omega| < 1$.
The absence of strong resonances means here the inequalities $\omega^2 \ne \{3/4, 8/9, 15/16\}$
to hold.

The system (\ref{8}) has an invariant symplectic four-dimensional
plane $q_3=q_4=$ $p_3=p_4=0$. The restriction of the system (\ref{8})
on this 4-plane gives the Hamiltonian system that corresponds to the stationary Swift-Hohenberg
equation on the spatial domain $\R$, and for $\alpha < 0$ its localized pulses match to
homoclinic orbits of the saddle-focus $O$ lying in this 4-plane. But for the system (\ref{8})
we are of interest in those homoclinic orbits of the equilibrium $O$ which do not lie in this
4-dimensional invariant plane. For such solution, if it exists, at least one of the coordinate
function $q_j(x), p_j(x),$ $j\ge 3,$ must not vanish identically.

Observe that the system (\ref{8}) is, in addition, reversible w.r.t. the linear
involution $L(q,p) = (L_q,L_p),$ $L_q = (q_1,-q_2,q_3,-q_4)$,
$L_p = (-p_1,p_2,-p_3,p_4).$ The fixed point set $Fix(L)$ of the involution $L$
is the 4-dimensional plane $q_2=q_4=$$p_1=p_3 =0.$ The equilibrium $O$ at the origin
of the Hamiltonian system (\ref{8}) is symmetric, i.e. $O\in Fix(L).$

In order to have homoclinic orbits of $O$, the equilibrium should possess eigenvalues with
positive and negative real parts. The linearization matrix at the equilibrium $O$ consists
of two independent $(4\times 4)$-blocks and its characteristic polynomial is the product
of two bi-quadratic polynomials
$$
[(\lambda^2+1)^2 -\alpha][(\lambda^2+1-\omega^2)^2 -\alpha].
$$
So, if $\alpha = 0$ and $\omega \in (1/2,1)$ the equilibrium has two double
pure imaginary eigenvalues $\pm i,$ $\pm i\sqrt{1-\omega^2}$, all of them are non
semi-simple ones, i.e. with $2\times 2$ Jordan blocks. For negative $\alpha$ roots
of the characteristic polynomial are two complex quartets
$$
\begin{array}{l}
\displaystyle{\pm\sqrt{\frac{\sqrt{1-\alpha}-1}{2}}\pm i\sqrt{\frac{\sqrt{1-\alpha}+1}{2}}},\\
\displaystyle{\pm\sqrt{\frac{\sqrt{(\omega^2-1)^2-\alpha}+\omega^2-1}{2}}\pm
i\sqrt{\frac{\sqrt{(\omega^2-1)^2-\alpha}+1-\omega^2}{2}}}.
\end{array}
$$
For negative $\alpha$ the squares of distances from the imaginary axis for the quartets
of eigenvalues in the complex plane $\mathbb C$ are ordered as follows
$$
\frac12[\sqrt{1-\alpha}-1] < \frac12[\sqrt{(\omega^2-1)^2-\alpha}+\omega^2 -1].
$$
So, the leading (corresponding to the eigenvalues closest to the imaginary axis)
stable direction (and unstable one, as well) of the equilibrium at
the origin is two-dimensional and coincides with the invariant plane corresponding
to the pair of eigenvalues with negative real parts of the first quartet. We expect
that homoclinic orbits of this equilibrium will approach to the
equilibrium along this direction (a general case).

Now we present the simulations which confirm the existence for small negative $\alpha$
of homoclinic orbits of $O$ not lying on the invariant 4-dimensional plane. We start
with some homoclinic orbit $\Gamma_0$ that does belong to the invariant 4-plane.
As is known theoretically \cite{GL,Dev,Knob,Sand} and by means of
simulations \cite{BGL,KLSh}, there are many homoclinic solutions of $O$ on
this invariant 4-plane. Having $\Gamma_0$ we examine the linearized system of
(\ref{8}) on this solution. This system is nonautonomous but asymptotically autonomous, it
splits into two linear four-dimensional subsystems due to the invariance of 4-plane, one
of which is invariant but another, the normal subsystem, contains the dependence on
$\Gamma_0$ in coefficients.

We study the normal variational 4-dimensional subsystem. Because the limit system has
the saddle-focus with the quartet of complex eigenvalues involving $\omega$ (see above),
the normal nonautonomous linear subsystem has a 3-dimensional invariant set of solutions which
asymptotically tend to zero solution, as $t\to \infty.$ Solutions of this set intersect any
cross-section $t=const > 0$ in the space $\R^4\times\R$ along a two-dimensional plane.
If we choose as $\Gamma_0$ a symmetric solution
w.r.t the involution $L$, then the normal linear subsystem will be reversible and for this
system there is another invariant set of solutions, obtained by the symmetry, they tend to zero
as $t\to -\infty.$ The intersection of corresponding two 2-planes at the cross-section $t=0$
is either transverse and therefore at only the origin, or they intersect along one-dimensional
subspace (the tangency), or they coincide. We seek for the intersection of these two-dimensional
planes calculating the related determinant composed
\begin{figure}[htp] 
\centering
\subfigure[]{%
\includegraphics[width=0.4\textwidth]{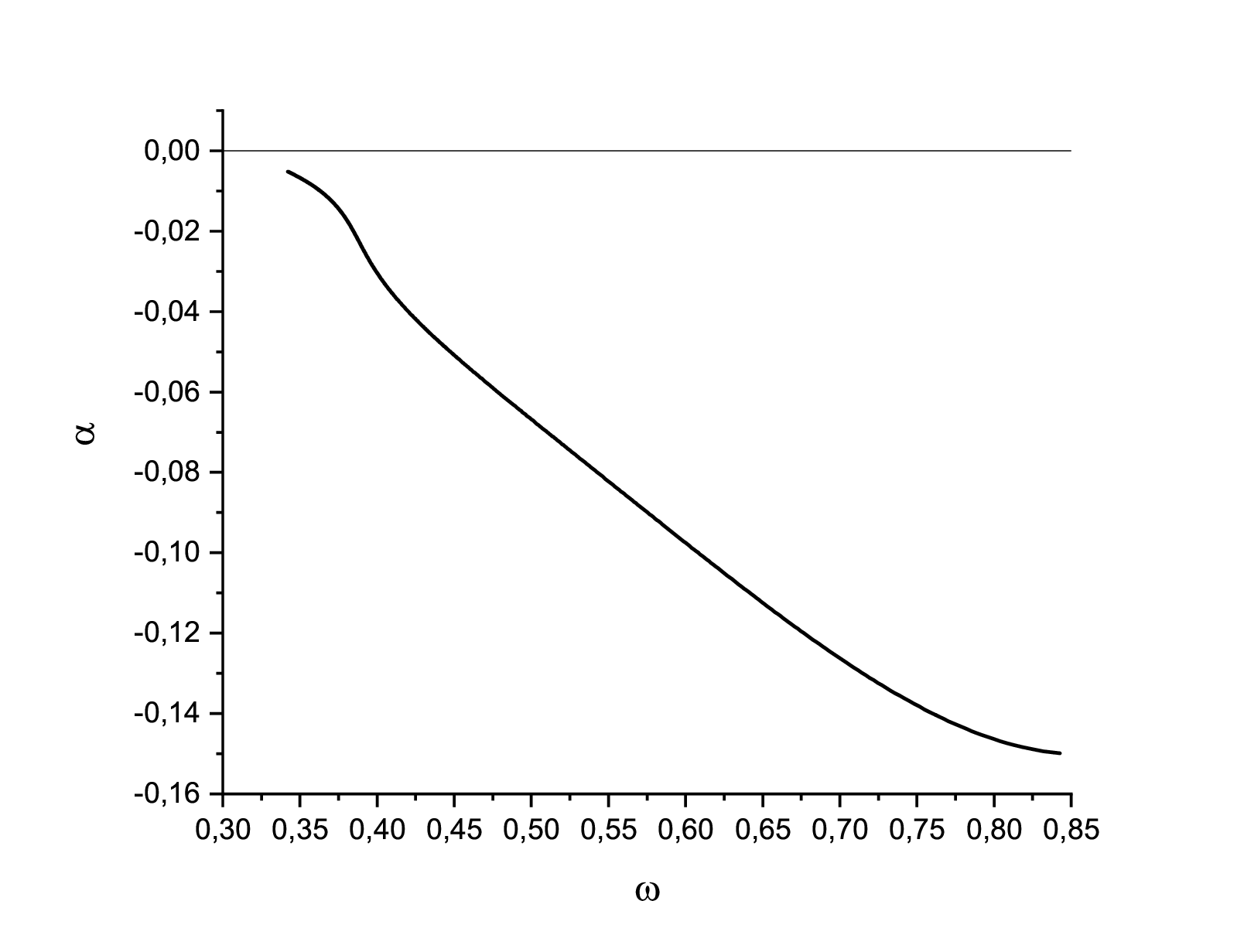}%
\label{fig:9:a}%
}\hfil
\subfigure[]{%
\includegraphics[width=0.4\textwidth]{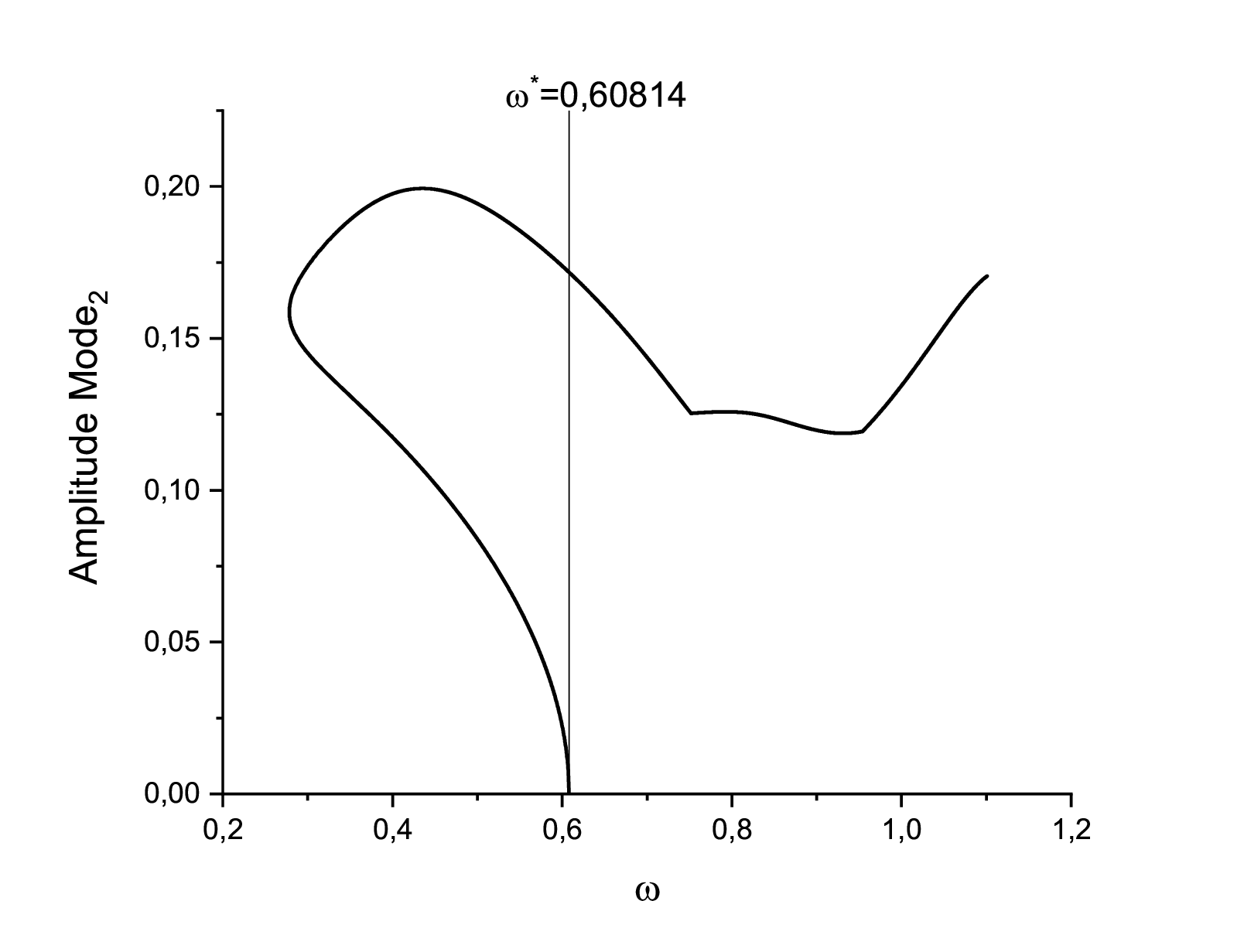}%
\label{fig:9:b}%
}
\caption{(a) The bifurcation curve of arising out-of-plane homoclinic solutions. (b) The dependence of bearing homoclinic solution on $\omega$.}
\label{fig:9}
\end{figure}

\begin{figure}[htp] 
\centering
\subfigure[]{%
\includegraphics[width=0.32\textwidth]{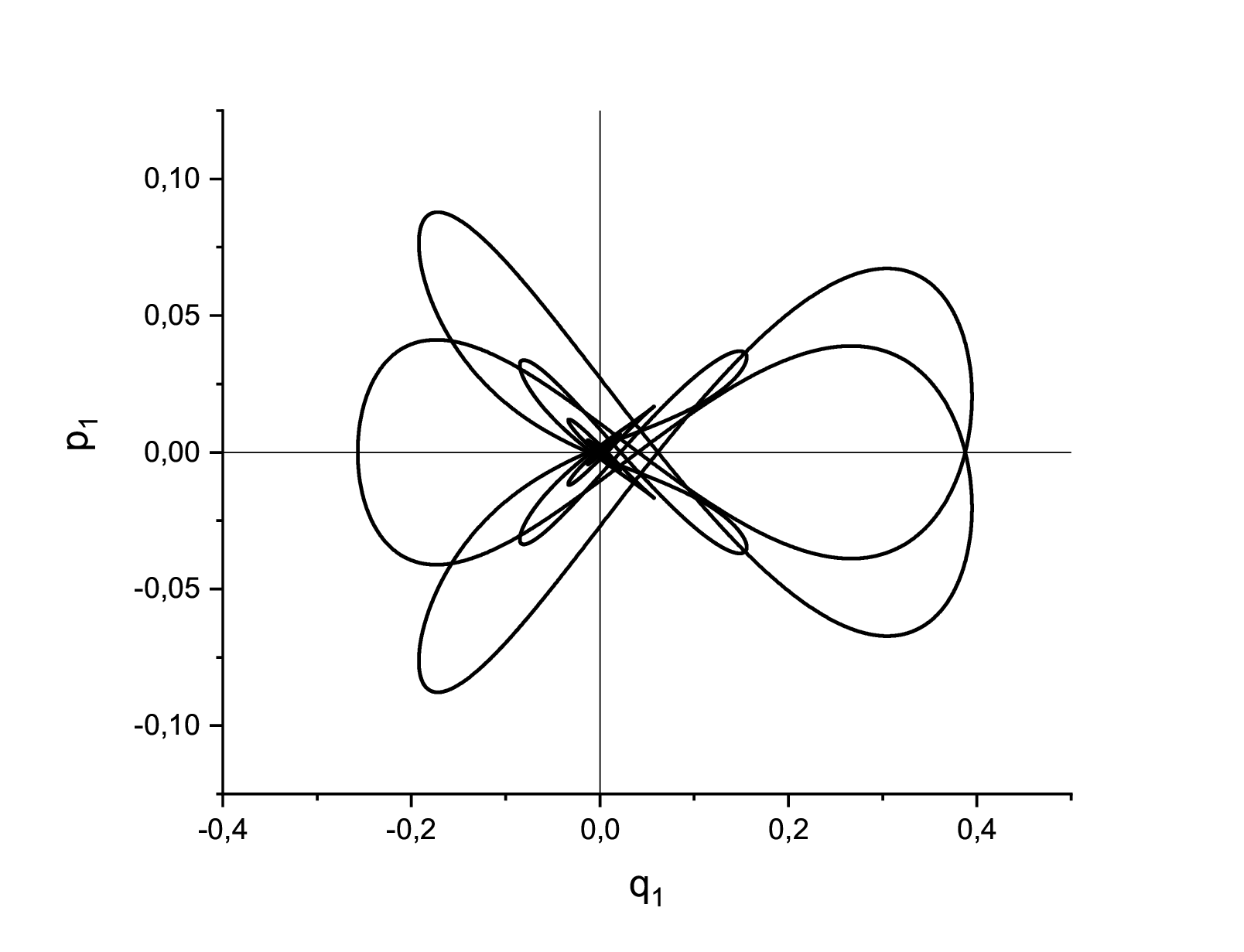}%
\label{fig:10:a}%
}\hfil
\subfigure[]{%
\includegraphics[width=0.32\textwidth]{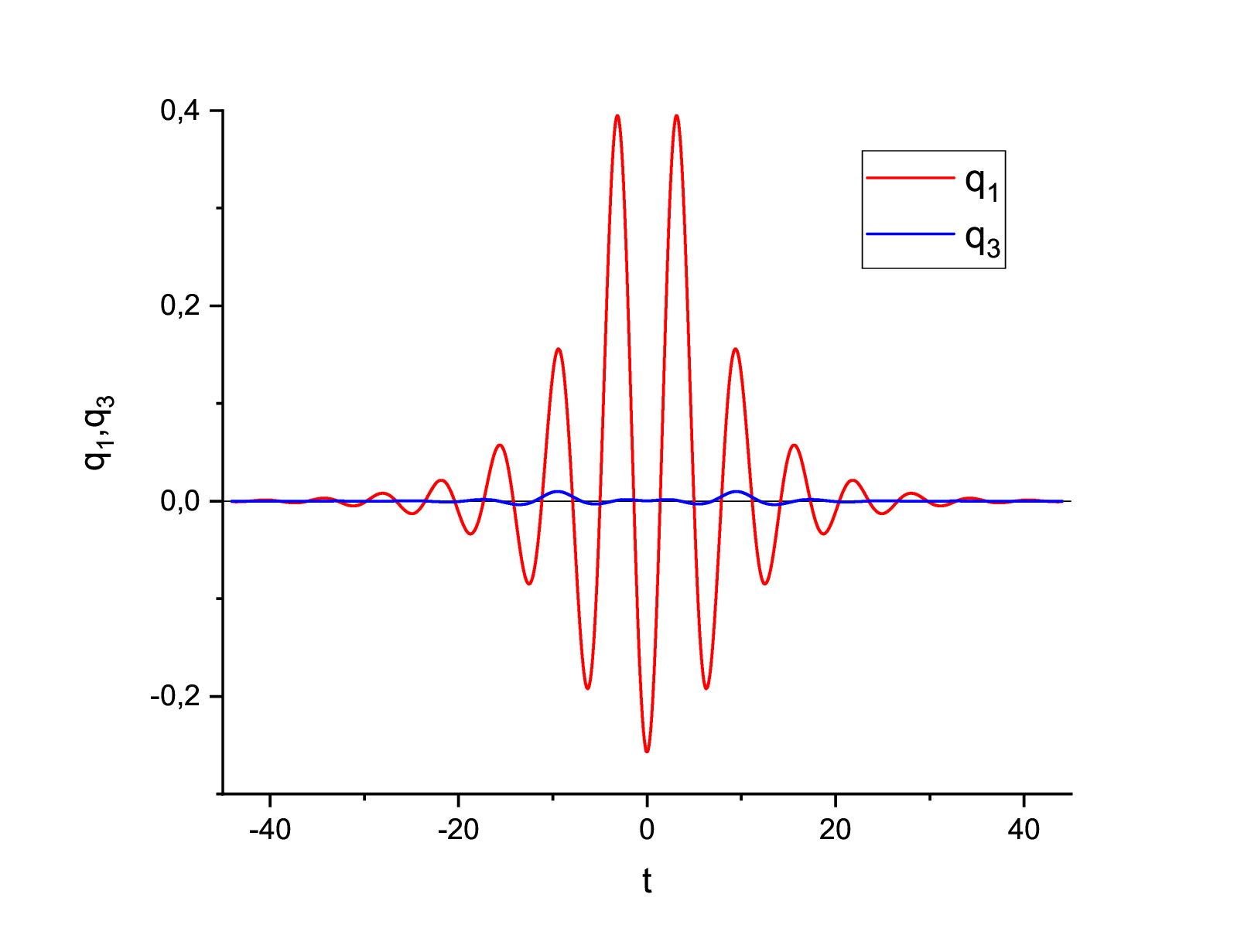}%
\label{fig:10:b}%
}
\subfigure[]{%
\includegraphics[width=0.32\textwidth]{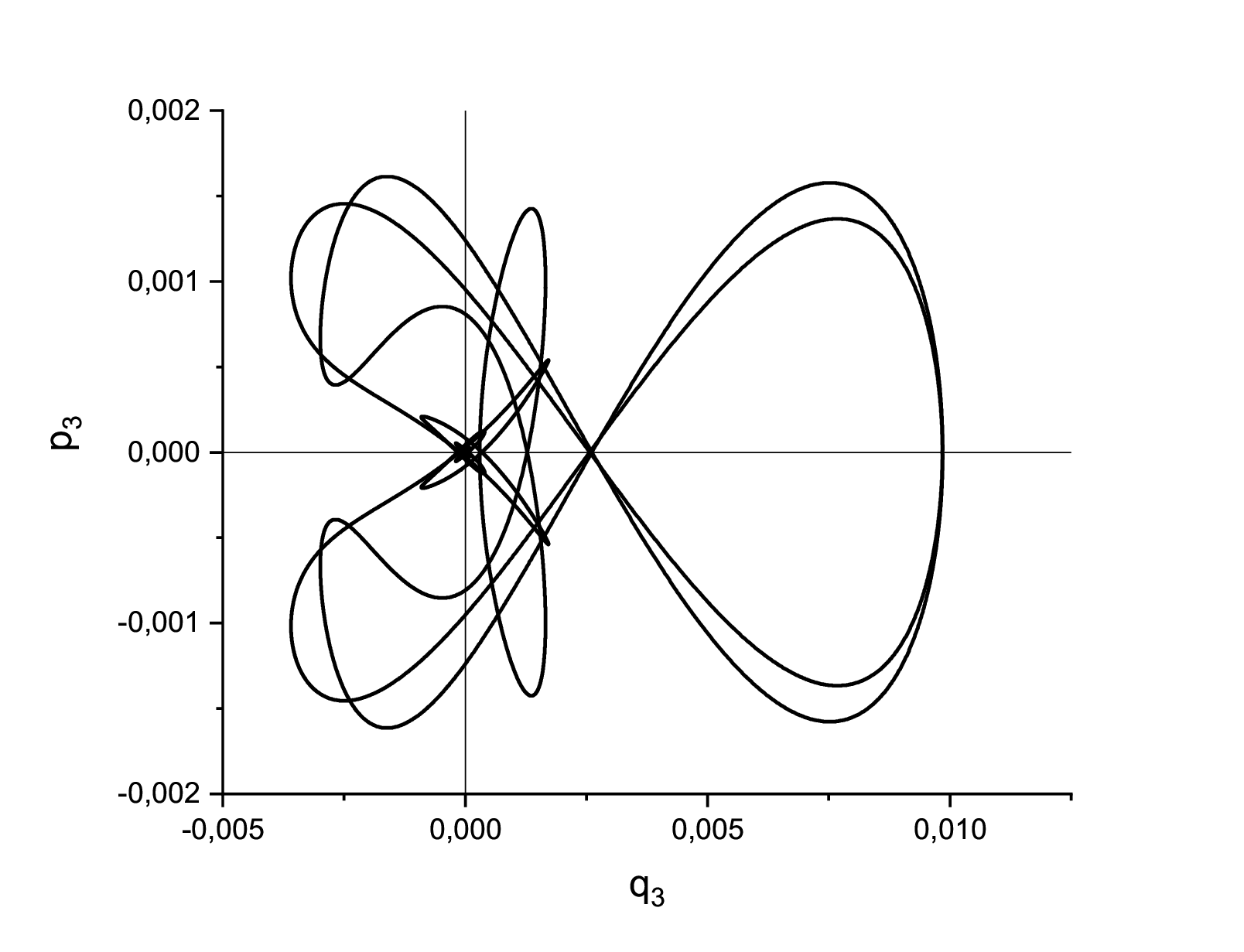}%
\label{fig:10:c}%
}\hfil
\caption{(a) $(p_1,q_1)$-projection, $\alpha=-0.1, \omega=0.60655$;
(b) unfoldings $q_1(t),q_3(t)$; (c) $(p_3,q_3)$-projection.}
\label{fig:10}
\end{figure}

\begin{figure}[htp] 
\centering
\subfigure[]{%
\includegraphics[width=0.32\textwidth]{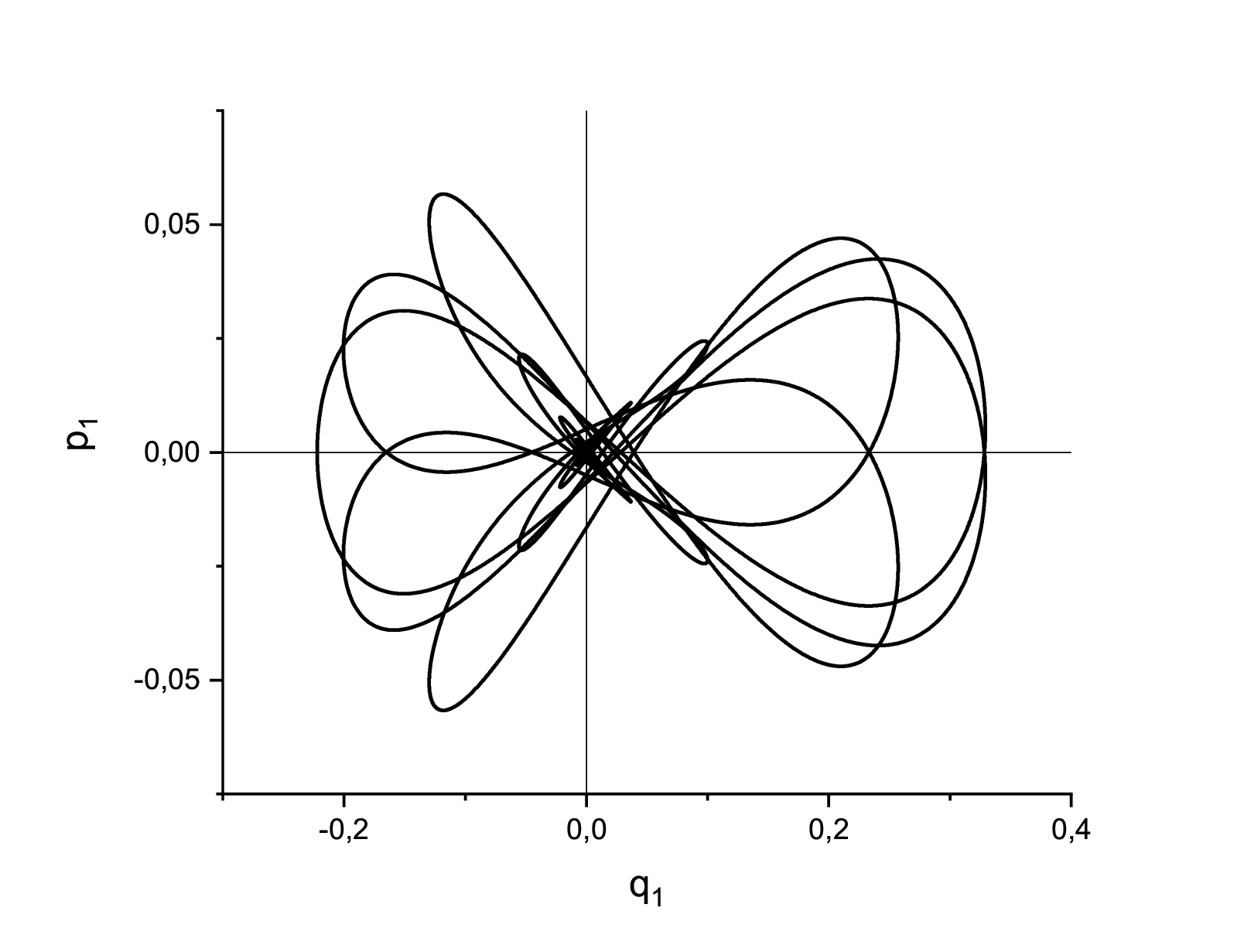}%
\label{fig:11:a}%
}\hfil
\subfigure[]{%
\includegraphics[width=0.32\textwidth]{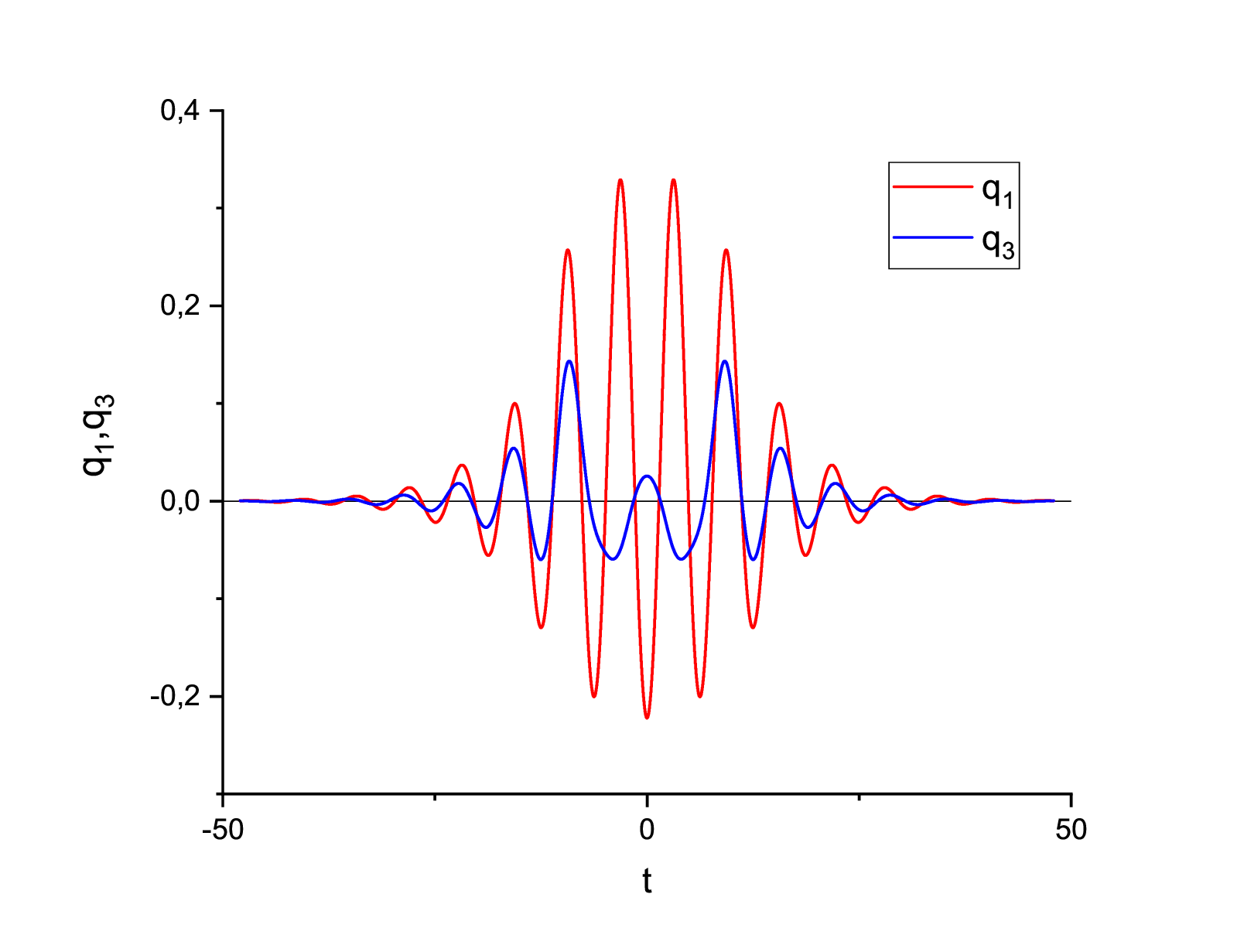}%
\label{fig:11:b}%
}
\subfigure[]{%
\includegraphics[width=0.32\textwidth]{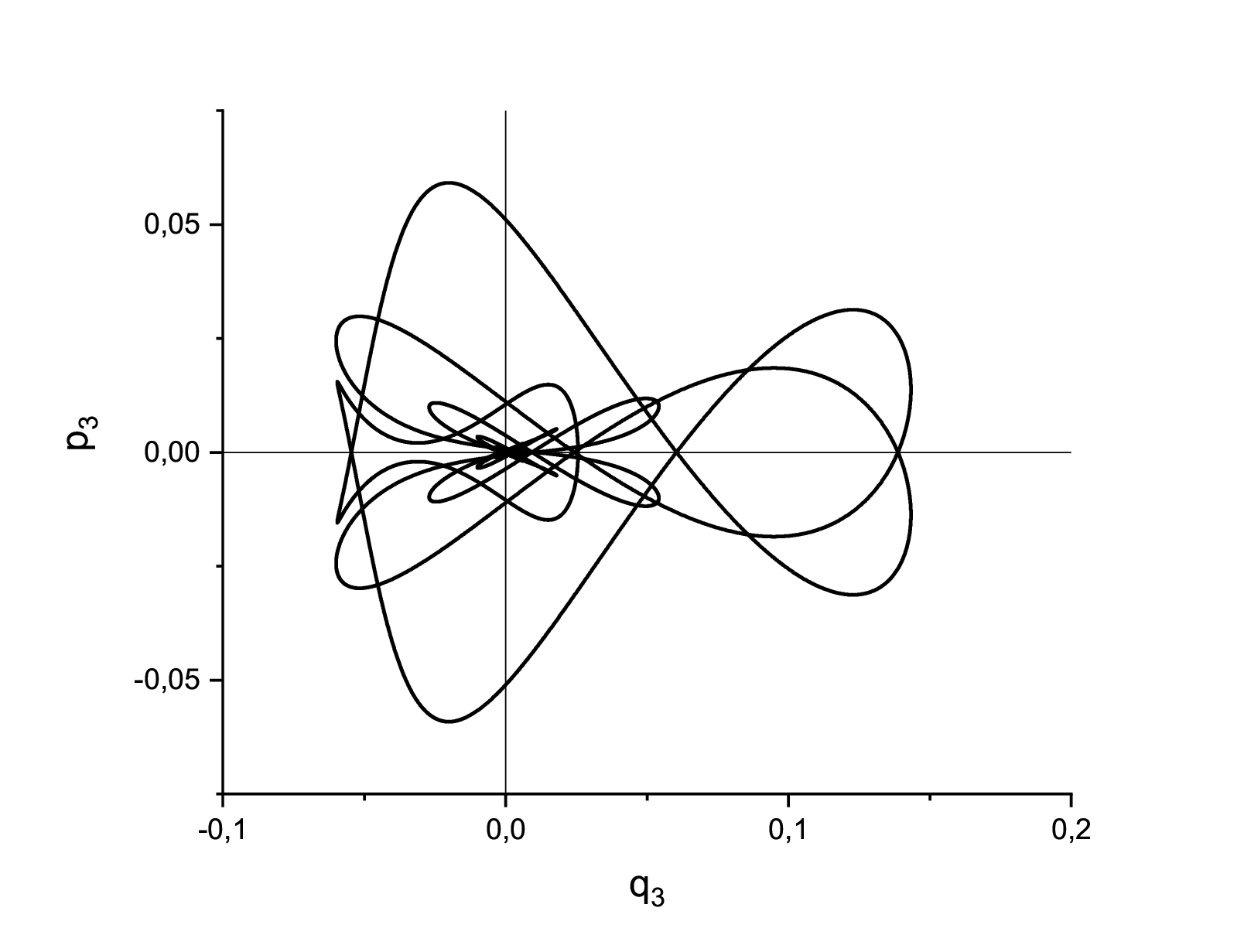}%
\label{fig:11:c}%
}\hfil
\caption{(a) $(p_1,q_1)$-projection, $\alpha=-0.1, \omega=0.29041$;
(b) unfoldings $q_1(t),q_3(t)$; (c) $(p_3,q_3)$-projection.}
\label{fig:11}
\end{figure}

from two independent vectors which span
the first plane and by symmetry have two other vectors in the second plane. The determinant
should have a simple zero at some value of the parameter $\omega.$ The related bifurcation
curve on the plane $(\alpha,\omega)$ is plotted on Fig.\ref{fig:9:a}.

The dependence of the norm of the born homoclinic out-of-plane solution of the
parameter $\omega$ is plotted on Fig.\ref{fig:9:b}.


After finding the intersection in the linear normal subsystem, that means the tangency
of stable and unstable manifolds of $O$ for the nonlinear system (\ref{8}), we search
numerically for the true homoclinic orbits outside of the invariant 4-plane. The results
are shown on the Fig.\ref{fig:10} and Fig.\ref{fig:11}.

\section{Conclusion}

The study we have performed shows that the double Hamiltonian Hopf Bifurcation is very reach if
one tries to investigate all orbit structures near the equilibrium for a generic 2-parameter
unfolding. We derived the normal form of the unfolding near the equilibrium, have found the
related invariants through which the normal form is expressed. After that we truncated the normal form
up to fourth order terms keeping in mind that this gives us the essential features of the orbit
structure. In this way we study some elements of the orbit structure of the reduced truncated system
in invariants as new variables. Also, we have proved that almost all systems in the unfolding
are nonintegrable and therefore even the dynamics of the truncated system is complicated.
To this end, we applied different tools developed in several papers, in particular the tool
related with separatrix splitting and the tool of differential Galois theory.
The nonintegrability is in contrast with the usual Hamiltonian Hopf bifurcation where the
truncated normal forms of any order is integrable. Also we discussed features which expect
to meet in the unfolding under the more detailed study. Here we made only first steps in this
way and many things have to be explored. We hope to proceed this study in the future work.

\section{Acknowledgement}

Authors thank D.V. Turaev and J.J. Morales-Ruiz for the valuable discussion and advices.
L.L.M. and N.E.K. acknowledge a support from the Russian Science Foundation (grant 22-11-00027).
This research (by R.M-S.) was in part supported by a grant from IPM (No. 1403700317).

\end{document}